\documentclass[leqno,a4paper,12pt]{article}

\usepackage{amsmath,amsthm}

\newtheorem{Thm}{\indent Theorem}[section]
\newtheorem{Prop}[Thm]{\indent Proposition}
\newtheorem{Lem}[Thm]{\indent Lemma}

\newtheorem{Cor}[Thm]{\indent Corollary}

\newtheorem{Var}[Thm]{\indent Variant}

\theoremstyle{definition}
\newtheorem{Def}[Thm]{\indent Definition}

\newtheorem{Rem}[Thm]{\indent Remark}
\newtheorem{Ex}[Thm]{\indent Example}

\newtheorem{Ass}[Thm]{\indent Assumption}

\def\qed{{\hskip0pt\unskip\unskip\nobreak\hfil\penalty50
          \hskip1em\hbox{}\nobreak\hfil
          {\bf q.e.d.}%
          \parfillskip=0pt\finalhyphendemerits=0
          \par}\medskip}

\newenvironment{Proof}
               {{\it Proof.}\quad}
               {\qed}

\newenvironment{Proofof}[1]
               {{\it Proof of #1.}\quad}
               {\qed}

\newcommand{\Prime}{\kern3\fontdimen1\font$'$\kern-7\fontdimen1\font}

\long\def\forget#1{}

\long\def\beginSIDEREMARK#1\endSIDEREMARK
    {{\par\bigskip\advance\leftskip by 2cm
                  \advance\rightskip by -2cm\noindent
      {\bf Our own side remark:} #1
      \par\bigskip\noindent}}
\long\def\beginFORGET#1\endFORGET{#1}
\long\def\beginFORGET#1\endFORGET{}

\def\?{\ ???\ \immediate\write16{}%
\immediate\write16{Warning: There was still a question mark . . . }%
\immediate\write16{}}

\usepackage{amsmath}
\usepackage{exscale}
\usepackage{amssymb}

\newcommand{\BC}{{\mathbb{C}}}

\newcommand{\BE}{{\mathbb{E}}}

\newcommand{\BG}{{\mathbb{G}}}

\newcommand{\BL}{{\mathbb{L}}}

\newcommand{\BP}{{\mathbb{P}}}
\newcommand{\BQ}{{\mathbb{Q}}}
\newcommand{\BR}{{\mathbb{R}}}

\newcommand{\BZ}{{\mathbb{Z}}}

\newcommand{\FH}{{\mathfrak{H}}}

\newcommand{\CD}{{\cal D}}

\newcommand{\CK}{{\cal K}}
\newcommand{\CL}{{\cal L}}

\newcommand{\CO}{{\cal O}}

\newcommand{\alg}{\mathop{\rm alg}\nolimits}
\newcommand{\Spec}{\mathop{{\bf Spec}}\nolimits}

\newcommand{\red}{\mathop{{\rm red}}\nolimits}
\newcommand{\reg}{\mathop{{\rm reg}}\nolimits}

\newcommand{\imm}{\mathop{{\rm im}}\nolimits}

\newcommand{\End}{\mathop{\rm End}\nolimits}
\newcommand{\Ext}{\mathop{\rm Ext}\nolimits}
\newcommand{\GL}{\mathop{\rm GL}\nolimits}
\newcommand{\Hom}{\mathop{\rm Hom}\nolimits}
\newcommand{\Pic}{\mathop{\rm Pic}\nolimits}

\newcommand{\loccit}{[loc.$\;$cit.]}

\def\halb{\frac{1}{2}}

\def\id{{\rm id}}

\newbox\mybox
\def\arrover#1{\mathrel{
       \setbox\mybox=\hbox spread 1.4em{\hfil$\scriptstyle#1$\hfil}
       \vbox{\offinterlineskip\copy\mybox
             \hbox to\wd\mybox{\rightarrowfill}}}}
\def\larrover#1{\mathrel{
       \setbox\mybox=\hbox spread 1.4em{\hfil$\scriptstyle#1$\hfil}
       \vbox{\offinterlineskip\copy\mybox
             \hbox to\wd\mybox{\leftarrowfill}}}}

\def\ontoover#1{\mathrel{
       \setbox\mybox=\hbox spread 1.4em{\hfil$\scriptstyle#1$\hfil}
       \vbox{\offinterlineskip\copy\mybox
             \hbox to\wd\mybox{\rightarrowfill\hskip-2.8mm
                               $\rightarrow$}}}}
\def\leftontoover#1{\mathrel{
       \setbox\mybox=\hbox spread 1.4em{\hfil$\scriptstyle#1$\hfil}
       \vbox{\offinterlineskip\copy\mybox
             \hbox to\wd\mybox{$\leftarrow$\hskip-2.8mm
                               \leftarrowfill}}}}
\def\longto{\longrightarrow}
\def\into{\hookrightarrow}
\def\onto{\ontoover{\ }}
\def\longonto{\ontoover{\ }}
\def\isoto{\arrover{\sim}}

\def\longinto{\lhook\joinrel\longrightarrow}

\usepackage[curve,matrix,arrow,cmtip]{xy}
\NoComputerModernTips

\def\myxymessage{\def\messagetext
   {Here an xy-pic diagram was omitted to speed up compilation . . . }
   \immediate\write16{\messagetext}
   \hbox{\bf \messagetext}}
\def\filxymatrix#1{\myxymessage}
\def\filxyarray#1{\myxymessage}
\newdir^{ (}{{}*!/-3pt/\dir^{(}}
\newdir_{ (}{{}*!/-3pt/\dir_{(}}
\newdir^{ )}{{}*!/+3pt/\dir^{)}}
\newdir_{ )}{{}*!/+3pt/\dir_{)}}

\def\rscript#1{\hbox to 0pt{$\scriptstyle#1$\hss}}

\let\oldbullet\bullet
\def\bullet{{\mathchoice{\oldbullet}%
                        {\oldbullet}%
                        {\scriptscriptstyle\oldbullet}%
                        {\oldbullet}}}

\newcommand{\argdot}{{\;\bullet\;}}%Punkt als Platzhalter fuer Argumente
\newcommand{\argast}{{\;\ast\;}}%Stern als Platzhalter fuer Argumente

\newcommand{\DM}{\mathop{DM^{eff}_-(k)}\nolimits}
\newcommand{\DeffgM}{\mathop{DM^{eff}_{gm}(k)}\nolimits}

\newcommand{\DeffQgM}{\mathop{DM^{eff}_{gm}(\BQ)}\nolimits}
\newcommand{\DgM}{\mathop{DM_{gm}(k)}\nolimits}
\newcommand{\ip}{\tilde{\imath}}
\newcommand{\jp}{\tilde{\jmath}}
\newcommand{\jss}{j_{!*}}
\newcommand{\KCE}{\mathop{{\rm KCE}}\nolimits}
\newcommand{\RC}{\mathop{{\bf R} C}\nolimits}
\newcommand{\ShN}{\mathop{Shv_{Nis}(SmCor(k))}\nolimits}

\newcommand{\ShDZ}{\mathop{Shv_{Zar}(D)}\nolimits}
\newcommand{\SmC}{\mathop{SmCor(k)}\nolimits}
\newcommand{\uC}{\mathop{\underline{C}}\nolimits}
\newcommand{\ap}{\mathop{\widetilde{a}} \nolimits}
\newcommand{\Xp}{\mathop{\widetilde{X}} \nolimits}
\newcommand{\Yp}{\mathop{\widetilde{Y}} \nolimits}
\newcommand{\bX}{\overline{X}}
\newcommand{\bY}{\overline{Y}}
\newcommand{\Xs}{\mathop{X^*} \nolimits}
\newcommand{\Ys}{\mathop{Y^*} \nolimits}
\newcommand{\wb}{\widetilde{\beta}}

\forget{

\newcommand{\Spec}{\mathop{{\rm Spec}}\nolimits}

}

\begin{document}

%%%%%%%%%%%%%%%%%%%%%%%%%%%%%%%%%%%%%%%%%%%%%%%%%%%%%%%%%%%%%%%%%%%%%%%
%
%  formatting

\hfuzz=3pt
\overfullrule=10pt                   % erzeugt schwarze Fehlerbalken

% The displayskip values were changed because \LaTeX does not react
% correctly to a \leqno: it should then use big skips, but doesn't.

\setlength{\abovedisplayskip}{6.0pt plus 3.0pt}
                               % preset 10.0pt plus 2.0pt minus 5.0pt
\setlength{\belowdisplayskip}{6.0pt plus 3.0pt}
                               % preset 10.0pt plus 2.0pt minus 5.0pt
\setlength{\abovedisplayshortskip}{6.0pt plus 3.0pt}
                               % preset 0.0pt plus 3.0pt
\setlength{\belowdisplayshortskip}{6.0pt plus 3.0pt}
                               % preset 6.0pt plus 3.0pt minus 3.0pt

\setlength{\baselineskip}{13.0pt}
                               % preset 12.0pt
\setlength{\lineskip}{0.0pt}
                               % preset 1.0pt
\setlength{\lineskiplimit}{0.0pt}
                               % preset 0.0pt

%%%%%%%%%%%%%%%%%%%%%%%%%%%%%%%%%%%%%%%%%%%%%%%%%%%%%%%%%%%%%%%%%%%%%%%
%
%  Title Page
%
%%%%%%%%%%%%%%%%%%%%%%%%%%%%%%%%%%%%%%%%%%%%%%%%%%%%%%%%%%%%%%%%%%%%%%%

\title{Pure motives, mixed motives and extensions of motives
associated to singular surfaces
\forget{
\footnotemark
\footnotetext{To appear in ....}
}
}
\author{\footnotesize by\\ \\
\mbox{\hskip-2cm
\begin{minipage}{6cm} \begin{center} \begin{tabular}{c}
J\"org Wildeshaus \\[0.2cm]
\footnotesize LAGA\\[-3pt]
\footnotesize UMR~7539\\[-3pt]
\footnotesize Institut Galil\'ee\\[-3pt]
\footnotesize Universit\'e Paris 13\\[-3pt]
\footnotesize Avenue Jean-Baptiste Cl\'ement\\[-3pt]
\footnotesize F-93430 Villetaneuse\\[-3pt]
\footnotesize France\\
{\footnotesize \tt wildesh@math.univ-paris13.fr}
\end{tabular} \end{center} \end{minipage}
\hskip-2cm}
%\\[2cm]
%{\bf Preliminary version --- not for distribution}\\[1cm]
}
% In the final version we might want to fix the date:
\date{May 12, 2011}
\maketitle
%\quad \\[-1.7cm]
\begin{abstract}
\noindent We first recall the construction of the Chow motive 
modelling intersection cohomology of a proper surface $\bX$, and
study its fundamental pro\-perties. Using Voevodsky's
category of effective geometrical motives, we then study the
motive of the exceptional divisor $D$ in a non-singular blow-up of $\bX$.
If all geometric irreducible components of $D$ are of genus
zero, then Voevodsky's formalism allows us to construct certain 
one-extensions of motives, as canonical 
sub-quotients of the motive with compact support 
of the smooth part of $\bX$. Speciali\-zing to Hilbert--Blumenthal
surfaces, we recover a motivic interpretation of a recent
construction of A.~Caspar.  \\

\noindent Keywords: intersection cohomology, intersection motives, 
K\"unneth filtration, motives \`a la Voevodsky, Kummer--Chern--Eisenstein
extensions.

%\noindent
%{\bf R\'esum\'e~:} RESUME.\\
\end{abstract}

%\vfill

\bigskip
\bigskip
\bigskip

\noindent {\footnotesize Math.\ Subj.\ Class.\ (2010) numbers: 14F42 
(11F41, 14C17, 14F43, 14G35, 14J99). }

\eject

\tableofcontents

\bigskip
%\vspace*{0.5cm}

%\newpage
%\include{Intro}
%%%%%%%%%%%%%%%%%%%%%%%%%%%%%%%%%%%%%%%%%%%%%%%%%%%%%%%%%%%%%%%%%%%%%%%
%
%  Introduction
%
%%%%%%%%%%%%%%%%%%%%%%%%%%%%%%%%%%%%%%%%%%%%%%%%%%%%%%%%%%%%%%%%%%%%%%%

\setcounter{section}{-1}
\section{Introduction}
\label{Intro}

%%%%%%%%%%%%%%%%%%%%%%%%%%%%%%%%

%%%%%%%%%%%%%%%%%%%%%%%%%%%%%%%%

The modest aim of this article is to construct 
non-trivial extensions in Voe\-vodsky's 
category of effective geometrical motives, by studying 
a very special and concrete geometric situation,
namely that of a singular proper surface. \\

This example illustrates a much more general
principle: varie\-ties $Y$ that are singular (or non-proper, for that matter),
can provide interesting extensions of motives. 
The cohomological theories of mixed sheaves suggest where to look
for these motives: the one should come from the open smooth
part $Y_{\reg}$ of $Y$ --- the \emph{intersection motive} of $Y$ --- 
the other
should be constructed out of the complement of $Y_{\reg}$ in 
(a compactification of) $Y$ --- 
the \emph{boundary motive} of $Y_{\reg}$.  
This principle (for which no originality
is claimed, since it has been part of the mathematical culture
for some time) will be discussed in more detail separately \cite{Wi}, 
in order to preserve the structure of the present article. 
It is intended as a research article with a large instructional
component. \\

The geometric object of interest is a proper surface $\bX$ over an 
arbitrary base field $k$. \\

The first three sections contain nothing fundamentally new,
except maybe for the systematic use of K\"unneth filtrations
(which are canonical) instead of K\"unneth decompositions
(which in general are not). Section~\ref{1} reviews a special
case of a result of Borho and MacPherson \cite{BoMp}, 
computing the intersection
cohomology of $\bX$ in terms of the cohomology of a 
desingularization $\Xp$. The
result, 
predicted by the Decomposition Theorem of \cite{BBD}, 
implies that the former is a direct
factor of the latter. More precisely (Theorem~\ref{1A}), its
complement is given by the second cohomology of the exceptional 
divisor $D$ of $\Xp$. 
As remarked already by de Cataldo and Migliorini
\cite{CM},
this fact can be interpreted motivically, which
allows one to construct the intersection motive
$h_{!*} (\bX)$ of $\bX$. This
is done in Section~\ref{2}. We get a canonical decomposition
\[
h(\Xp) = h_{!*} (\bX) \oplus \bigoplus_m h^2(D_m)
\]
in the category of Chow motives over $k$. 
Recall that this category is
pseudo-Abelian. The above decomposition
should be considered as remarkable: to construct a sub-motive of $h(\Xp)$ 
does not \emph{a priori}
necessitate the \emph{identification}, but only the 
\emph{existence} of a complement. In our
situation, the complement \emph{is} canonical, thanks to the very 
special geometrical situation.
This point is reflected by the rather subtle functoriality properties 
of $h_{!*} (\bX)$
(Proposition~\ref{2E}): viewed as a sub-motive of $h(\Xp)$, 
it is respected by pull-backs,
viewed as a quotient, it is respected by push-forwards 
under dominant morphisms of surfaces.
Section~\ref{3} is devoted to the existence and the study of 
the K\"unneth filtration of $h_{!*}(\bX)$. The main ingredient is of course 
Murre's construction of K\"unneth projectors for the motive 
$h(\Xp)$ \cite{Mr}. 
Theorem~\ref{3C} shows how to adapt these to our construction. \\

As suggested by one of the fundamental properties of intersection
cohomology \cite{BBD}, the intersection motive of $\bX$ satisfies the Hard 
Lefschetz Theorem for ample line bundles on $\bX$. We prove this result
(Theorem~\ref{4A}) in Section~\ref{4}. In fact, we give a slight
generalization (Variant~\ref{4A'}), which will turn out to be useful
for the setting we shall study in the last section. \\ 

Section~\ref{5} is concerned with the motive of the
boundary $D$ of the desingularization $\Xp$ of $\bX$. This boundary
being singular in general, the right language for the study of 
its motive is given by Voevodsky's triangulated category of
effective geometrical motives \cite{VSF}. The section starts with a review
of the definition of this category, and of its relation to Chow
motives. It is then easy to define motivic analogues of $H^0$ and $H^2$
of $D$, and to see that they are Chow motives. The most interesting
part is the motivic analogue of the part of degree one 
$H^1$, which will be seen as a canonical
sub-quotient of the motive of $D$. \\

In Section~\ref{6}, we unite what was said before, and give our main
result (Theorem~\ref{Main}). Assuming that all geometric
irreducible components 
of $D$ are of genus zero, we construct a
one-extension of the degree two-part of the intersection
motive of $\bX$ by the degree one-part of the motive of $D$. 
We have no difficulty to admit that this statement was greatly
inspired by the main result
of a recent article of Caspar \cite{Cs}. It thus appeared appropriate to
conclude this article by a discussion of his result. This
is done in Section~\ref{7}, where we show that in the geometric
setting considered in \loccit , Theorem~\ref{Main} yields a
motivic interpretation of Caspar's construction. \\

Part of this work was done while I was enjoying a 
\emph{cong{\'e} pour recherches ou conversions th{\'e}matiques},
granted by the \emph{Universit{\'e} Paris~13}, and during a visit to the \emph{Centre de Recerca
Matem\`atica} at Bellaterra--Barcelona. I am grateful to both institutions. I also wish to thank J.\ Ayoub, J.-B.~Bost,
J.I.\ Burgos, M.A.A.\ de Cataldo, F.\ D\'eglise,
B.\ Kahn, K.\ K\"unnemann and F.\ Lemma for useful comments and discussions. \\

{\bf Notations and convention}: $k$ denotes a fixed base field, 
and $CH$ stands for the tensor
product with $\BQ$ of the Chow group. The $\BQ$-linear
category of Chow motives over
$k$ is denoted by $CHM(k)_{\BQ}$.
Our standard reference for Chow motives is Scholl's survey
article \cite{Sch}.

%%% Local Variables:
%%% mode: latex
%%% TeX-master: "imot"
%%% End:

\bigskip

%\include{Sec1}
%%%%%%%%%%%%%%%%%%%%%%%%%%%%%%%%%%%%%%%%%%%%%%%%%%%%%%%%%%%%%%%%%%%%%%%
%
%  Section 1
%
%%%%%%%%%%%%%%%%%%%%%%%%%%%%%%%%%%%%%%%%%%%%%%%%%%%%%%%%%%%%%%%%%%%%%%%

\section{Intersection cohomology of surfaces}
\label{1}

%%%%%%%%%%%%%%%%%%%%%%%%%%%%%%%%

%%%%%%%%%%%%%%%%%%%%%%%%%%%%%%%%

In order to motivate the construction of the intersection motive, 
to be given in the next section, we
shall recall the computation
of the \emph{intersection cohomology} of a complex surface. \\

Thus, throughout this section, our base field $k$ will be equal to $\BC$. We consider the
following situation:
\[
\vcenter{\xymatrix@R-10pt{
        X \ar@{^{ (}->}[r]^-{j} &
        \Xs \ar@{<-^{ )}}[r]^{i} &
        Z \; .
\\}}
\]
The morphism $i$ is a closed immersion of a sub-scheme $Z$, with complement $j$. The scheme $X^*$
is a surface over $\BC$, all of whose singularities are contained in $Z$. Thus, the surface
$X$ is smooth. \\

Our aim is to identify the intersection cohomology groups $H^n_{!*} (\Xs(\BC),\BQ)$. Note that
since $X$ is smooth, the complex $\BQ_X [2]$ consisting of the constant local system $\BQ$, placed
in degree $-2$, can be viewed as a \emph{perverse sheaf} 
(for the middle perversity) on $X(\BC)$ 
\cite[Sect.~2.2.1]{BBD}. Hence its \emph{intermediate extension} $\jss \BQ_X [2]$
\cite[(2.2.3.1)]{BBD} is defined as a perverse sheaf on $X^*(\BC)$. By definition,
\[
H^n_{!*} (\Xs(\BC),\BQ) = H^{n-2} (X^*(\BC),\jss \BQ_X [2]) \; , \; \forall \, n \in \BZ \; .
\]
In order to identify $H^n_{!*} (\Xs(\BC),\BQ)$, note first that the normalization of $\Xs$ is
finite over $\Xs$, and the direct image under finite morphisms is exact for the perverse
$t$-structure \cite[Cor.~2.2.6~(i)]{BBD}. Therefore, intersection cohomology is invariant under
passage to the normalization. In the sequel, we therefore assume $\Xs$ to be normal. In
particular, its
singularities are isolated. \\

Next, note that if $\Xs$ is smooth, then the complex $\jss \BQ_X [2]$ equals $\BQ_{\Xs} [2]$.
Transitivity of $\jss$ \cite[(2.1.7.1)]{BBD} shows that we may enlarge $X$, and hence assume that
the closed sub-scheme $Z$ is finite. \\

Choose a resolution of singularities. More precisely, consider in addition the following diagram,
assumed to be cartesian:
\[
\vcenter{\xymatrix@R-10pt{
        X \ar@{^{ (}->}[r]^-{\jp} \ar@{=}[d] &
        {\Xp} \ar@{<-^{ )}}[r]^{\ip} \ar[d]_\pi &
        D \ar[d]^\pi \\
        X \ar@{^{ (}->}[r]^-{j} &
        \Xs \ar@{<-^{ )}}[r]^{i} &
        Z
\\}}
\]
The morphism $\pi$ is assumed proper (and birational) and the surface 
$\Xp$, smooth. We then have the following special case of
\cite[Thm.~1.7]{BoMp}.

\begin{Thm} \label{1A}
(i) For $n \ne 2$,
\[
H^n_{!*} (\Xs(\BC),\BQ) = H^n (\Xp(\BC),\BQ) \; .
\]
\noindent (ii) The group $H^2_{!*} (\Xs(\BC),\BQ)$ is a direct factor of $H^2 (\Xp(\BC),\BQ)$, with
a \emph{cano\-ni\-cal} complement. As a sub-group, this complement is given by the map
\[
{\ip}_*: H^2_{D(\BC)}(\Xp(\BC),\BQ) \longto H^2 (\Xp(\BC),\BQ)
\]
from cohomology with support in $D(\BC)$; this map is injective. As a quotient, the complement is
given by the restriction
\[
{\ip}^*: H^2 (\Xp(\BC),\BQ) \longto H^2(D(\BC),\BQ) \; ;
\]
this map is surjective.
\end{Thm}

Note that this result is compatible with further blow-up of $\Xp$ in points belonging to $D$. \\

Let us construct the maps between $H^n_{!*}
(\Xs(\BC),\BQ)$
and $H^n (\Xp(\BC),\BQ)$ leading to the above identifications. 
Consider the total direct image $\pi_* \BQ_{\Xp}$~; following the convention used in \cite{BBD}, we
drop the letter ``$R$'' from our notation.

\begin{Lem} \label{1B}
The complex $\pi_* \BQ_{\Xp}[2]$ is a perverse sheaf on $X^*$.
\end{Lem}

\begin{Proof}
Let $P$ be a point (of $Z$) over which $\pi$ is not an isomorphism, and denote by $i_P$ its
inclusion into $\Xs$. By definition \cite[D{\'e}f.~2.1.2]{BBD}, we need to check that (a)~the higher
inverse images $H^n i_P^* \pi_* \BQ_{\Xp}$ vanish for $n>2$, (b)~the higher exceptional inverse
images $H^n i^!_P \pi_* \BQ_{\Xp}$ vanish for $n<2$.

(a) By proper base change, the group in question equals $H^n (\pi^{-1}(P),\BQ)$. Since
$\pi^{-1}(P)$ is of dimension at most one, there is no cohomology above degree two.

(b) The surface $\Xp$ is smooth. Duality and proper base change imply that the group in question is
abstractly isomorphic to the dual of $H^{4-n} (\pi^{-1}(P),\BQ)$. This group vanishes if $4-n$ is
strictly larger than two.
\end{Proof}

For $a \in \BZ$, denote by $\tau_{\le a}$ the functor associating to a complex the $a$-th step of
its canonical filtration (with respect to the classical $t$-structure). Recall that $\jss \BQ_X
[2]$ equals $\tau_{\le -1} (j_* \BQ_X[2])$ \cite[Prop.~2.1.11]{BBD}. We now see how to relate it
to $\pi_* \BQ_{\Xp}[2]$: apply $\tau_{\le -1} \circ \pi_*$ to the exact triangle
\[
{\ip}_* {\ip}^! \BQ_{\Xp} [2] \longto \BQ_{\Xp} [2] \longto {\jp}_* \BQ_X [2] \longto 
{\ip}_* {\ip}^! \BQ_{\Xp} [3] \; .
\]
We get an exact triangle
\[
i_* \tau_{\le -1} (i^! \pi_* \BQ_{\Xp}[2] ) \longto \tau_{\le -1} ( \pi_* \BQ_{\Xp}[2] ) \longto 
\jss \BQ_X [2] \longto i_* \tau_{\le -1} (i^! \pi_* \BQ_{\Xp}[2]) [1] \; .
\]
But according to Lemma~\ref{1B}, 
\[
\tau_{\le -1} (i^! \pi_* \BQ_{\Xp}[2]) = 0 \; .
\]
We thus get the following.

\begin{Lem} \label{1C}
There is a canonical isomorphism
\[
\tau_{\le -1} (\pi_* \BQ_{\Xp}[2]) \isoto \jss \BQ_X [2] 
\]
of perverse sheaves on $X^*$.
\end{Lem}

\begin{Proofof}{Theorem~\ref{1A}}
By Lemma~\ref{1C}, there is a canonical exact triangle
\[
\jss \BQ_X [2] \longto \pi_* \BQ_{\Xp}[2] \longto (\tau_{\ge 2} \pi_* \BQ_{\Xp}) [2] \longto
(\jss \BQ_X [2])[1] \; .
\]
It implies that
\[
H^n_{!*} (\Xs(\BC),\BQ) = H^n (\Xp(\BC),\BQ) 
\]
for $n = 0, 1$, and that the sequence
\[
0 \longto H^2_{!*} (\Xs(\BC),\BQ) \longto H^2 (\Xp(\BC),\BQ) \stackrel{{\ip}^*}{\longto}
H^2 (D(\BC),\BQ)
\]
is exact. Duality implies that
\[
H^n_{!*} (\Xs(\BC),\BQ) = H^n (\Xp(\BC),\BQ) 
\]
for $n = 3,4$, too. Therefore, the sequence
\[
0 \longto H^2_{!*} (\Xs(\BC),\BQ) \longto H^2 (\Xp(\BC),\BQ) \stackrel{{\ip}^*}{\longto}
H^2 (D(\BC),\BQ) \longto 0
\]
is exact. Hence so is the dual exact sequence
\[
0 \longto H^2_{D(\BC)}(\Xp(\BC),\BQ)\stackrel{{\ip}_*}{\longto}  H^2 (\Xp(\BC),\BQ) \longto
H^2_{!*} (\Xs(\BC),\BQ) \longto 0 \; .
\]
\end{Proofof}

\begin{Rem} \label{1D}
The analogue of Theorem~\ref{1A} holds for $\ell$-adic cohomology, 
and when $k$ is a
finite field of characteristic unequal to $\ell$. The proof is exactly the 
same. Note that by
Abhyankar's result on resolution of singularities in dimension two \cite[Theorem]{L2}, $X^*$ can
be desingularized for \emph{any} base field $k$. In addition 
(see the discussion in \cite[pp.~191--194]{L1}), by further blowing
up possible singularities of (the components of) the pre-image $D$ of $Z$, 
it can be assumed to be
a divisor with
normal crossings, whose irreducible components are smooth.
This discussion also shows that the system of such resolutions is filtering.
\end{Rem}

\begin{Rem}
Theorem~\ref{1A}~(ii) implies that the composition
\[
{\ip}^* {\ip}_* : H^2_{D(\BC)}(\Xp(\BC),\BQ) \longto H^2(D(\BC),\BQ)
\]
is an isomorphism. 
\end{Rem}

\bigskip

%\include{Sec2}
%%%%%%%%%%%%%%%%%%%%%%%%%%%%%%%%%%%%%%%%%%%%%%%%%%%%%%%%%%%%%%%%%%%%%%%
%
%  Section 2
%
%%%%%%%%%%%%%%%%%%%%%%%%%%%%%%%%%%%%%%%%%%%%%%%%%%%%%%%%%%%%%%%%%%%%%%%

\section{Construction of the intersection motive}
\label{2}

%%%%%%%%%%%%%%%%%%%%%%%%%%%%%%%%

%%%%%%%%%%%%%%%%%%%%%%%%%%%%%%%%

Fix a base field $k$, and assume given a proper surface $\bX$ over $k$. The aim of this section is
to recall the construction of 
the \emph{Chow motive} modelling intersection cohomo\-logy of $\bX$,
and to study its functoriality properties. The discussion
preceding Theorem~\ref{1A} showed 
that intersection cohomology is invariant under passage to the
normalization $\Xs$ of $\bX$; the same should thus be expected from 
the motive we intend to construct.
\footnotemark \footnotetext{ This principle also explains why the problem of constructing the
intersection motive of a proper curve $\overline{C}$ is not very interesting: 
the intersection motive of
$\overline{C}$ is equal to the motive of the normalization $C^*$ of $\overline{C}$ (which is
smooth and projective).} Fix
\[
\vcenter{\xymatrix@R-10pt{
        X \ar@{^{ (}->}[r] &
        \Xs \ar@{<-^{ )}}[r]^{i} &
        Z
\\}}
\]
where $i$ is a closed immersion of a finite sub-scheme $Z$, with smooth complement $X$. Choose a
resolution of singularities. More precisely, consider in addition the following diagram, assumed
to be cartesian:
\[
\vcenter{\xymatrix@R-10pt{
        X \ar@{^{ (}->}[r] \ar@{=}[d] &
        \Xp \ar@{<-^{ )}}[r]^{\ip} \ar[d]_\pi &
        D \ar[d]^\pi \\
        X \ar@{^{ (}->}[r] &
        \Xs \ar@{<-^{ )}}[r]^{i} &
        Z
\\}}
\]
where $\pi$ is proper (and birational), $\Xp$ is smooth (and proper), and $D$ is a divisor with
normal crossings, whose irreducible components $D_m$ are smooth (and proper). \\
\begin{Rem} \label{2a}
Note that $\Xp$, as a smooth and proper surface, is projective:
Zariski proved this result
for algebraically closed base fields in \cite[p.~54]{Z},
and \cite[Cor.~7.7]{SGA} allows to descend to arbitrary base fields.
\end{Rem}
\medskip

Theorem~\ref{1A} suggests how to construct the intersection motive;
in particular, it should be a canonical direct
complement of $\oplus_m h^2(D_m)$ in $h(\Xp)$. 
Recall \cite[Sect.~1.13]{Sch} that the $h^2(D_m)$ are canonically defined as quotient
objects of the motives $h(D_m)$. Hence there is a canonical morphism
\[
{\ip}^*: h(\Xp) \longto \bigoplus_m h(D_m) \longonto \bigoplus_m h^2(D_m)
\]
of Chow motives. Similarly \cite[Sect.~1.11]{Sch}, 
there is a canonical morphism
\[
{\ip}_*: \bigoplus_m h^0(D_m)(-1) \longinto \bigoplus_m h(D_m)(-1) \longto h(\Xp) \; .
\]
Here, the twist by $(-1)$ denotes the tensor product with the Lefschetz motive $\BL=h^2(\BP^1)$. The following is a special case of
\cite[Sect.~2.5]{CM}.

\begin{Thm} \label{2B}
(i) The composition $\alpha := {\ip}^*{\ip}_*$ is an isomorphism of Chow motives. \\[0.1cm]
(ii) The composition $p:= {\ip}_*\alpha^{-1}{\ip}^*$ is an idempotent on $h(\Xp)$.
Hence so is the difference $\id_{\Xp}-p$. \\[0.1cm]
(iii) The image $\imm p$ is canonically isomorphic to $\oplus_m h^2(D_m)$.
\end{Thm}

\begin{Proof}
(ii) and (iii) are formal consequences of (i). 
The formula ``$\phi_* \phi^* = \deg \phi$'' for
finite morphisms $\phi$ \cite[Sect.~1.10]{Sch} 
shows that we may prove our claim after a finite
extension of our ground field $k$. In particular, we may assume that all components $D_m$ are
geometrically irreducible, with field of constants equal to $k$. We then have canonical
isomorphisms $h^0(D_m) \cong h(\Spec k)$ and $h^2(D_m) \cong \BL$. Denote by $i_m$ the closed
immersion of $D_m$ into $\Xp$. The map $\alpha$ in question equals
\[
\bigoplus_{m,n} \;  i_m^* i_{n,*}: \bigoplus_n h^0(D_n)(-1) \longto \bigoplus_m h^2(D_m) \; .
\]
For each pair $(m,n)$, the composition $i_m^* i_{n,*}$ is an endomorphism of $\BL$. Now the degree
map induces an isomorphism
\[
\End (\BL) = CH^0 (\Spec k) \isoto \BQ \; .
\]
We leave it to the reader to show that under this isomorphism, the endomorphism $i_m^* i_{n,*}$ is
mapped to the intersection number $D_n \cdot D_m$. Our claim follows from the non-degeneracy of
the intersection pairing on the components of $D$ \cite[p.~6]{M}.
\end{Proof}

Following \cite[p.~158]{CM}, we propose the following definition.

\begin{Def} \label{2C}
The \emph{intersection motive} of $\bX$ is defined as
\[
h_{!*} (\bX) := (\Xp,\id_{\Xp}-p,0) \in CHM(k)_{\BQ} \; .
\]
\end{Def}

Here, we follow the standard notation for Chow motives (see e.g.\ 
\cite[Sect.~1.4]{Sch}). Idempotents on
Chow motives admit an image; by definition, the image of the idempotent $\id_{\Xp}-p$ on the Chow
motive $(\Xp,\id_{\Xp},0)=h(\Xp)$ is $(\Xp,\id_{\Xp}-p,0)=h_{!*} (\bX)$. Note that by definition, we
have
the equality $h_{!*} (\bX)=h_{!*} (\Xs)$. \\

Theorem~\ref{2B} shows that there is a canonical decomposition
\[
h(\Xp) = h_{!*} (\bX) \oplus \bigoplus_m h^2(D_m)
\]
in $CHM(k)_{\BQ}$.
By Theorem~\ref{1A} and Remark~\ref{1D}, the Betti, resp.\ $\ell$-adic realization of the
intersection motive (for the base fields for which this realization exists) coincides with
intersection cohomology of $\bX$ (and of $\Xs$).

\begin{Prop} \label{2D}
As before, denote by $\Xs$ the normalization of $\bX$. The definition of $h_{!*} (\bX)$ is
independent of the choices of the finite sub-scheme $Z$ containing the singularities, and of
the desingularization $\Xp$ of $\Xs$.
\end{Prop}

This statement is going to be proved together with the functoriality pro\-perties of the
intersection motive, whose formulation we prepare now. Consider a dominant morphism $f:\bX \to
\bY$ of proper surfaces over $k$. By the universal property of the normalization $\Ys$ of $\bY$,
it induces a morphism, still denoted $f$, between $\Xs$ and $\Ys$. It is generically finite. Hence
we can find a finite closed subscheme $W$ of $\Ys$ containing the singularities, and such that the
pre-image under $f$ of $Y := \Ys - W$ is dense, and smooth. The closed sub-scheme $f^{-1}(W)$ of
$X$ contains the singularities of $\Xs$. We thus can find a morphism $F$ of desingularizations of
$\Xs$ and $\Ys$ of the type considered before:
\[
\quad\quad\quad\quad (F) \quad\quad\quad\quad\quad\quad\quad\quad
\vcenter{\xymatrix@R-10pt{
        \Xp \ar@{<-^{ )}}[r]^{i_D} \ar[d]_F &
        D \ar[d]^F \\
        \Yp \ar@{<-^{ )}}[r]^{i_C} &
        C
\\}}
\quad\quad\quad\quad\quad\quad\quad\quad\quad\quad\quad\quad
\]
This means that $\Xp$ and $\Yp$ are smooth, and $D$ and $C$ are divisors with normal crossings,
whose irreducible components $D_m$ resp.\ $C_n$
are smooth, and lying over finite closed sub-schemes of $\Xs$ and
$\Ys$, respectively. Choose and fix such a dia\-gram. Note that if the original morphism $f:\bX \to
\bY$ is finite, then the diagram $(F)$ can be chosen to be cartesian: first, choose $\Yp$
and define $f_1 : \bX_1 \to \Yp$ as the base change
$\bX \times_{\bY} \Yp$ of $\Yp$ via the morphism $f$. The latter being finite,
the irreducible components of $f^{-1}_1(C)$ lie over finite closed sub-schemes of $\bX$.
The surface $\Xp$ is then obtained by further blowing up $\bX_1$. 

\begin{Prop} \label{2E}
(i) The pull-back $F^*: h(\Yp) \to h(\Xp)$ maps the sub-object $h_{!*}(\bY)$ of $h(\Yp)$ to the
sub-object $h_{!*}(\bX)$ of $h(\Xp)$. \\[0.1cm]
(ii) The push-forward $F_*: h(\Xp) \to h(\Yp)$ maps the quotient $h_{!*}(\bX)$ of $h(\Xp)$ to the
quotient $h_{!*}(\bY)$ of $h(\Yp)$. \\[0.1cm]
(iii) The composition $F_*F^*: h_{!*}(\bY) \to h_{!*}(\bY)$ equals multiplication with the degree
of
$f$. \\[0.1cm]
(iv) If $f$ is finite, and if the diagram $(F)$ is chosen to be cartesian, then both $F^*$ and
$F_*$ respect the decompositions
\[
h(\Yp) = h_{!*} (\bY) \oplus \bigoplus_n h^2(C_n)
\]
and
\[
h(\Xp) = h_{!*} (\bX) \oplus \bigoplus_m h^2(D_m)
\]
of $h(\Yp)$ and of $h(\Xp)$, respectively.
\end{Prop}

\begin{Proof}
By definition, there are (split) exact sequences
\[
0 \longto h_{!*} (\bX) \longto h(\Xp) \stackrel{i_D^*}{\longto} 
\bigoplus_m h^2(D_m) \longto 0
\]
and
\[
0 \longto \bigoplus_m h^0(D_m)(-1) \stackrel{i_{D,*}}{\longto} h(\Xp) \longto h_{!*} (\bX) \longto
0 \; ;
\]
similarly for $\Yp$ and $C$. Obviously, the first sequence is contravariant, and the second is
covariant. This proves parts (i) and (ii). Part (iii) follows from this, and from the
corresponding formula for $F_*F^*$ on the motive of $\Yp$ 
\cite[Sect.~1.10]{Sch}; note that the degree of
$F$ equals the one of $f$. If $(F)$ is cartesian, then the above sequences are both co- and
contravariant thanks to the base change formulae $F_*i_D^* = i_C^*F_*$ and $F^*i_{C,*} =
i_{D,*}F^*$. This proves part (iv).
\end{Proof}

\begin{Proofof}{Proposition~\ref{2D}}
First, let us show that for a fixed choice of $Z$, the definition of $h_{!*} (\bX)$ is independent
of the choice of the desingularization $\Xp$ of $\Xs$. Using that the system of such
desingularizations is filtering, we reduce ourselves to the situation considered in
Proposition~\ref{2E}, with $f=\id$. We thus have a cartesian diagram
\[
\vcenter{\xymatrix@R-10pt{
        \Xp \ar@{<-^{ )}}[r]^{i_D} \ar[d]_F &
        D \ar[d]^F \\
        \Xp' \ar@{<-^{ )}}[r]^{i_C} &
        C
\\}}
\]
Let us denote by $h_{!*} (\bX)$ and $h_{!*}' (\bX)$ the two intersection motives formed with
respect to $\Xp$ and $\Xp'$, respectively. We want to show that $F^*: h_{!*}' (\bX) \to
h_{!*}(\bX)$ is an isomorphism. The scheme $\Xp'$ is normal, and the morphism $F$ is proper. By
the valuative criterion of properness, the locus of points of $\Xp'$ where $F^{-1}$ is not defined
is of dimension zero. 
By \cite[Prop.~V.5.3]{Ha}, $\Xp$ dominates the blow-up of $\Xp'$ in the points $P_1,\ldots,P_r$ where
$F$ is not an isomorphism. This blow-up lies between $\Xp$ and $\Xp'$, and satisfies the same
conditions on desingularizations. Repeating this argument and using the fact that 
the number of irreducible components of the fibres $F^{-1}(P_i)$ is
finite, one sees that this process stops at some stage; $F$ is therefore the composition of
blow-ups in points. By
induction, we may assume that $F$ equals the blow-up of $\Xp'$ in
one point $P$. The exceptional divisor $E :=
F^{-1}(P)$ is a projective bundle (of rank one) over $P$. It is also one of the irreducible
components $D_m$ of $D$; 
in fact, the morphism $F$ induces a bijection between the components of $D$
other than $E$ and the components $C_n$ of $C$. 
Denote by $i_E$ the closed immersion of $E$ into $\Xp$.
By Manin's computation of the motive of a blow-up \cite[Thm.~2.8]{Sch}, the sequence
\[
0 \longto h(\Xp') \stackrel{F^*}{\longto} h(\Xp) \stackrel{i_E^*}{\longto} h^2(E) \longto 0
\]
is (split) exact. But obviously, so is
\[
0 \longto \bigoplus_n h^2(C_n) \stackrel{F^*}{\longto} 
\bigoplus_m h^2(D_m) \stackrel{i_E^*}{\longto} h^2(E) \longto 0 \; .
\]
Hence $F^*$ maps the kernel $h_{!*}' (\bX)$ of $i_C^*$ isomorphically to the kernel $h_{!*}
(\bX)$ of $i_D^*$.

In the same way, one shows that enlarging $Z$ by adding non-singular points of $X^*$ does not
change the value of $h_{!*} (\bX)$.
\end{Proofof}

Recall the definition of the \emph{dual} of a Chow motive 
\cite[Sect.~1.15]{Sch}. For example, for any
desingularization $\Xp$ of $X^*$, the dual of $(\Xp,\id_{\Xp},0)=h(\Xp)$ is given by
$(\Xp,\id_{\Xp},2)=h(\Xp)(2)$.

\begin{Prop} \label{2F}
The dual of the intersection motive $h_{!*} (\bX)$ is canonically isomorphic to $h_{!*} (\bX)(2)$.
\end{Prop}

\begin{Proof}
By definition, the dual of $(\Xp,\id_{\Xp}-p,0)$ equals $(\Xp,{ }^t(\id_{\Xp}-p),2)$, where ${ }^t$
denotes the transposition of cycles in $\Xp \times \Xp$. But $p$ is symmetric: in fact, ${}^t({\ip}^*)
= {\ip}_*$, and ${}^t({\ip}_*) = {\ip}^*$.

One checks as in the proof of Proposition~\ref{2D} that this identification of $h_{!*} (\bX)^*$
with $h_{!*} (\bX)(2)$ does not depend on the choice of $\Xp$.
\end{Proof}

%%% Local Variables:
%%% mode: latex
%%% TeX-master: "imot"
%%% End:

\bigskip

%\include{Sec3}
%%%%%%%%%%%%%%%%%%%%%%%%%%%%%%%%%%%%%%%%%%%%%%%%%%%%%%%%%%%%%%%%%%%%%%%
%
%  Section 3
%
%%%%%%%%%%%%%%%%%%%%%%%%%%%%%%%%%%%%%%%%%%%%%%%%%%%%%%%%%%%%%%%%%%%%%%%

\section{The K\"unneth filtration of the intersection motive}
\label{3}

%%%%%%%%%%%%%%%%%%%%%%%%%%%%%%%%

%%%%%%%%%%%%%%%%%%%%%%%%%%%%%%%%

We continue to consider the situation of Section~\ref{2}. Thus, $\bX$ is a proper surface over the
base field $k$ with normalization $\Xs$, and we fix
\[
\vcenter{\xymatrix@R-10pt{
        X \ar@{^{ (}->}[r] &
        \Xs \ar@{<-^{ )}}[r]^{i} &
        Z
\\}}
\]
where $i$ is a closed immersion of a finite sub-scheme $Z$, 
with smooth complement $X$. In
addition, we consider the following cartesian diagram:
\[
\vcenter{\xymatrix@R-10pt{
        X \ar@{^{ (}->}[r] \ar@{=}[d] &
        \Xp \ar@{<-^{ )}}[r]^{\ip} \ar[d]_\pi &
        D \ar[d]^\pi \\
        X \ar@{^{ (}->}[r] &
        \Xs \ar@{<-^{ )}}[r]^{i} &
        Z
\\}}
\]
where $\pi$ is proper, $\Xp$ is smooth and proper (hence projective), 
and $D$ is a divisor with normal crossings,
whose
irreducible components $D_m$ are smooth. The aim of this section is to 
recall Murre's
construction of \emph{K\"unneth decompositions} of the motive of $\Xp$ 
\cite{Mr}, following
Scholl's presentation \cite[Chap.~4]{Sch},
and to study the resulting filtration on the intersection motive. \\

Thus, fix (i)~a hyperplane section $C \subset \Xp$ that is a
smooth curve (observe that $C$ might only be defined over a finite
extension $k'$ of $k$). As explained in \cite[Sect.~4.3]{Sch}, 
the embedding of $C$ into $\Xp$ induces an
isogeny $P \to J$ from the Picard variety to the Albanese variety of $\Xp$. 
This isogeny is
actually independent of the choice of 
the smooth curve $C$ representing the fixed very ample class in 
$CH^1(\Xp)$ (and a non-zero multiple of the isogeny 
is defined over $k$). 
Fix (ii)~an isogeny $\beta: J \to P$ such that the
composition of the two isogenies equals multiplication by $n > 0$. Finally, fix (iii)~a $0$-cycle
$T$ of degree one on $C$. Then by \cite[Thm.~3.9]{Sch}, $\beta$ corresponds to a symmetric cycle
class
\[
\wb \in CH^1 (\Xp \times \Xp)
\]
satisfying the condition $p_{\Xp,*} (\wb \cdot [\Xp \times T]) = 0 \in CH^1 (\Xp)$, where $p_{\Xp}$ is
the first projection from the product $\Xp \times \Xp$ to $\Xp$. \\

One then defines \cite[Sect.~4.3]{Sch} 
projectors $\pi_0 := [T \times \Xp]$ and $\pi_4 := { }^t \pi_0 =
[\Xp \times T]$, as well as $p_1 := \frac{1}{n} \wb \cdot [C \times \Xp]$ and $p_3 := { }^t p_1$.
All orthogonality relations are satisfied, including $p_3 p_1 = 0$, 
except that $p_1 p_3$ is not
necessarily equal to zero. This is why a modification is necessary: one puts $\pi_1 := p_1 - \halb
p_1 p_3$ and $\pi_3 := { }^t \pi_1 = p_3 - \halb p_1 p_3$. \footnotemark \footnotetext{ This
differs from Murre's original solution \cite[Rem.~6.5]{Mr}, where one takes $p_1 - p_1 p_3$ and
$p_3$ instead of $\pi_1$ and $\pi_3$.} This, together with $\pi_2 := \id_{\Xp} - \pi_0 - \pi_1 -
\pi_3 - \pi_4$, gives a full auto-dual set of orthogonal projectors. We thus get a K\"unneth
decomposition of $h(\Xp)$
(first over $k'$, then by pushing down, over $k$):
\[
h(\Xp) = {}'h^0(\Xp) \oplus {}'h^1(\Xp) \oplus {}'h^2(\Xp) 
\oplus {}'h^3(\Xp) \oplus {}'h^4(\Xp) \; ,
\]
with
\[
{}'h^n(\Xp) := (\Xp, \pi_n, 0) \subset (\Xp, \id_{\Xp}, 0) = h(\Xp) \; ,
\quad 0 \le n \le 4 \; .
\]
\begin{Def} \label{3a}
(a)~The \emph{K\"unneth
filtration of $h(\Xp)$} is the ascending filtration of $h(\Xp)$ 
by sub-motives induced by a K\"unneth
decomposition of $h(\Xp)$:
\[
0 \subset h^0(\Xp) \subset h^{\le 1}(\Xp) \subset h^{\le 2}(\Xp) 
\subset h^{\le 3}(\Xp) \subset h^{\le 4}(\Xp) = h(\Xp) \; ,
\]
where we set $h^{\le r} (\Xp) := \oplus_{n=0}^r {}'h^n(\Xp)$, 
$r \le 4$. \\[0.1cm]
(b)~The $n$-th \emph{K\"unneth component of $h(\Xp)$}, $0 \le n \le 4$,
is the sub-quotient of $h(\Xp)$ defined by
\[
h^n(\Xp) := h^{\le n}(\Xp) / h^{\le n-1}(\Xp) \; .
\]
\end{Def}

\begin{Rem}
The sub-objects $h^{\le n}(\Xp)$ are direct factors of $h(\Xp)$, hence
the sub-quotients $h^n(\Xp)$ exist.
Similarly, one may define the quotients 
\[
h^{\ge r} (\Xp) := h(\Xp) / h^{\le {r-1}} (\Xp) 
\] 
of $h(\Xp)$.
\end{Rem}

Note that a number of choices is involved in the construction of the
projectors $\pi_0,\ldots,\pi_4$: mainly, a very ample line bundle $\CL$
on $\Xp$,
and a $0$-cycle on a smooth curve in the divisor class corresponding to $\CL$. 
The following is the content of
\cite[Thm.~14.3.10~i)]{KMrP}. 

\begin{Prop} \label{3b}
The K\"unneth filtration of $h(\Xp)$ is independent of the choices
made in the construction of the K\"unneth decomposition.
\end{Prop}

\begin{Rem} \label{3c}
(a)~In particular, the K\"unneth components $h^n(\Xp)$ are ca\-nonically
defined sub-quotients of $h(\Xp)$. \\[0.1cm]
(b)~\emph{A posteriori}, one may define the notion of K\"unneth
decomposition of $h(\Xp)$ as being a decomposition splitting the K\"unneth
filtration. Such decompositions include the ones obtained by  
Murre's construction, but there could be others.
\end{Rem}

Our aim (see Theorem~\ref{3C})
is to deduce from the K\"unneth filtration of $h(\Xp)$ 
a filtration of the intersection motive 
$h_{!*}(\bX) \subset h(\Xp)$:
\[
0 \subset h^0_{!*}(\bX) \subset h^{\le 1}_{!*}(\bX) \subset h^{\le 2}_{!*}(\bX) 
\subset h^{\le 3}_{!*}(\bX) \subset h^{\le 4}_{!*}(\bX) = h_{!*}(\bX) \; .
\]
The idea is of course to take the ``induced'' filtration.
But since we are working in a category which is only pseudo-Abelian,
we need to proceed with some care.   
Recall the quotient
$\oplus_m h^2 (D_m)$ and the sub-object $\oplus_m h^0 (D_m)$ of $\oplus_m h (D_m)$.

\begin{Prop} \label{3B}
The K\"unneth filtration of $h(\Xp)$ satisfies the following
conditions.
\begin{enumerate}
\item[(1)] Duality $h(\Xp)^{\vee} \isoto h(\Xp)(2)$ induces isomorphisms
\[
h^{\le r}(\Xp)^{\vee} \isoto h^{\ge 4-r}(\Xp)(2) \; .
\]
\item[(2)] The composition of morphisms
\[
h^{\le 1}(\Xp) \longinto h(\Xp) \stackrel{{{\ip}}^*}{\longto} \bigoplus_m h (D_m) \longonto \bigoplus_m
h^2 (D_m)
\]
equals zero.
\end{enumerate}
\end{Prop}

\begin{Proof}
The K\"unneth filtration satisfies (1) since the decompositions obtained by 
Murre's construction are auto-dual: 
${}'h^n(\Xp)^{\vee} \cong {}'h^{4-n}(\Xp)(2)$  under the duality 
$h(\Xp)^{\vee}
\cong h(\Xp)(2)$. 

By \cite[Prop.~5.8]{J}, condition~(2) is a
consequence of Murre's Conjecture~B \cite[Sect.~1.4]{Mr2} on the triviality of the action of the
$\ell$-th K\"unneth projector on $CH^j (Y)$, for $\ell > 2j$. Here, $Y$ equals the product of $\Xp$
and $D_m$, $j=2$, and $\ell = 5, 6$. Note that for products of a surface and a curve, the
conjecture is known to hold (see \cite[Lemma~8.3.2]{Mr3} for the case $j=2$).

But since the argument proving (2) is rather explicit, we may just as well give it for the
convenience of the reader. We need to compute the composition of correspondences
\[
h(\Xp) \stackrel{\pi_n}{\longto} h(\Xp) \stackrel{{{\ip}}^*}{\longto} \bigoplus_m h (D_m)
\stackrel{pr}{\longonto} \bigoplus_m h^2 (D_m) \; ,
\]
for $n = 0, 1$. The composition is zero if and only if it 
is zero after base change to a finite field
extension. Hence we may assume that all $D_m$ are geometrically irreducible, with field of
constants $k$. Then the $h^2 (D_m)$ equal $\BL$, and the composition $pr \circ {{\ip}}^*$ corresponds
to the cycle class
\[
([D_m])_m \in \bigoplus_m CH^1(\Xp)
\]
on $\coprod_m \Xp \times \Spec k$. By definition of the composition of correspondences, we then
find
\[
pr \circ {{\ip}}^* \circ \pi =  \bigl( p_{\Xp,*}(\pi \cdot [\Xp \times D_m]) \bigr)_m \in \bigoplus_m
CH^1(\Xp) \; ,
\]
for any $\pi \in CH^2(\Xp \times \Xp)$. Here as before, $p_{\Xp}$ is the first projection from the
product $\Xp \times \Xp$ to $\Xp$. Let us fix $m$. We need to show that for $n = 0, 1$, the cycle
class
\[
p_{\Xp,*}(\pi_n \cdot [\Xp \times D_m]) \in CH^1(\Xp)
\]
is zero. For $n = 0$, this is easy: the intersection
\[
\pi_0 \cdot [\Xp \times D_m] = [T \times \Xp] \cdot [\Xp \times D_m] = [T \times D_m]
\]
has one-dimensional fibres under $p_{\Xp}$. Therefore, its push-forward under $p_{\Xp}$ is zero.

For $n=1$, observe first that by definition of $\pi_1$, and by associativity of composition of
correspondences, it suffices to show that
\[
p_{\Xp,*}(p_1 \cdot [\Xp \times D_m]) = 0 \; .
\]
By definition, the intersection $p_1 \cdot [\Xp \times D_m]$ is a non-zero multiple of
\[
\wb \cdot [C \times \Xp] \cdot [\Xp \times D_m] \; .
\]
By the projection formula, the image under $p_{\Xp,*}$ of this cycle equals the image under the
push-forward $CH^0(C) \to CH^1(\Xp)$ of
\[
pr_{1,*} ( \wb_C \cdot [C \times D_m] ) \; ,
\]
where $\wb_C$ denotes the pull-back of $\wb$ to $C \times \Xp$, and $pr_1$ the projection from $C
\times \Xp$ to $C$. Denote by $pr_2$ the projection from this product to $\Xp$. Now symmetry of $\wb$
and the condition $p_{\Xp,*} (\wb \cdot [\Xp \times T]) = 0$ imply that
\[
pr_{2,*}(\wb_C \times [T \times \Xp]) = 0 \in CH^1(\Xp) \; .
\]
It follows that
\[
pr_{2,*}(\wb_C \times [T \times D_m]) = 0 \in CH^1(D_m) \; ,
\]
where we denote by the same symbol $pr_2$ the projection from
$C \times D_m$ to $D_m$.
In particular, the degree $a$ of this $0$-cycle is zero. But since $T$ is of degree one, we have
\[
pr_{1,*} ( \wb_C \cdot [C \times D_m] ) = a [C] \in CH^0(C) \; .
\]
\end{Proof}

Given that duality $h(D_m)^{\vee} \isoto h(D_m)(1)$ induces an isomorphism
\[
h^0(D_m)^{\vee} \isoto h^2(D_m)(1) \; ,
\]
it is easy to see that the morphism ${\ip}_*$ dual to the one from condition~(2)
\[
\bigoplus_m h^0 (D_m) \longinto \bigoplus_m h (D_m) \stackrel{{\ip}_*}{\longto} h(\Xp)(1) \longonto
h^{\ge 3}(\Xp)(1)
\]
is zero, i.e., the map ${\ip}_*: \oplus_m h^0 (D_m) \to h(\Xp)(1)$ factors through the sub-motive
$h^{\le 2}(\Xp)(1)$. On the other hand, by condition~(2), the inverse image ${{\ip}}^*: h(\Xp) \to
\oplus_m h^2 (D_m)$ factors through the quotient motive $h^{\ge 2}(\Xp)$. It follows that the
composition
\[
\alpha = {{\ip}}^* {\ip}_* : \bigoplus_m h^0 (D_m)(-1) \longto \bigoplus_m h^2 (D_m)
\]
considered in Section~\ref{2} factors naturally through $h^2(\Xp)$. 
By Theorem~\ref{2B}~(i), the morphism $\alpha$ is an isomorphism. 

\begin{Def} \label{3d}
Define the motive $h^2_{!*}(\bX)$ as the kernel of 
\[
{\ip}_* \alpha^{-1} {{\ip}}^* : h^2(\Xp) \longto h^2(\Xp) \; .
\]
\end{Def}

Note that ${\ip}_* \alpha^{-1} {{\ip}}^*$ 
is an idempotent on $h^2(\Xp)$; it therefore
admits a kernel. Its image is of course canonically isomorphic 
(via ${{\ip}}^*$) to
$\oplus_m h^2 (D_m)$. Dually, the image of the projector 
$\id_{h^2(\Xp)} - {\ip}_* \alpha^{-1} {{\ip}}^*$ is $h^2_{!*}(\bX)$.
Its kernel is canonically isomorphic 
(via ${{\ip}}_*$) to
$\oplus_m h^0 (D_m)(-1)$.

\begin{Rem}
In \cite[Sect.~14.2.2]{KMrP}, 
the \emph{transcendental part} $t^2(\Xp)$ of the motive of the surface $\Xp$ is
defined, as a complement
in $h^2(\Xp)$ of the algebraic, i.e., ``N\'eron--Severi''-part $h^2(\Xp)_{\alg}$. It follows that
under the projection from $h^2(\Xp)$, the transcendental part $t^2(\Xp)$
maps monomorphically to $h^2_{!*}(\bX)$.  
\end{Rem}

By condition~(2) from
Proposition~\ref{3B}, the projector $p = {\ip}_*\alpha^{-1}{\ip}^*$ 
on $h(\Xp)$ used to define $h_{!*}(\bX)$ 
gives rise to compatible factorizations
\[
p^{\ge r} := {\ip}_*\alpha^{-1}{\ip}^* : 
h^{\ge r}(\Xp) \longto h^{\ge r}(\Xp) \; , \; r \le 2
\]
and 
\[
p^{\le r} := {\ip}_*\alpha^{-1}{\ip}^* : 
h^{\le r}(\Xp) \longto h^{\le r}(\Xp) \; , \; r \ge 2 \; ,
\]
all of which are again idempotent. Consequently, we get (split) exact
sequences of motives
\[
0 \longto h^{\le 1}(\Xp) \longto \ker (p^{\le 2})
\longto h^2_{!*}(\bX) \longto 0 \; ,
\]
\[
0 \longto \ker (p^{\le 2})  \longto \ker (p^{\le 3})
\longto h^3(\Xp) \longto 0 
\]
etc.

\begin{Thm} \label{3C}
(i)~The K\"unneth filtration
of $h(\Xp)$ 
\[
0 \subset h^0(\Xp) \subset h^{\le 1}(\Xp) \subset h^{\le 2}(\Xp) 
\subset h^{\le 3}(\Xp) \subset h^{\le 4}(\Xp) = h(\Xp) 
\]
induces a filtration of the intersection motive 
$h_{!*}(\bX)$
\[
0 \subset h^0_{!*}(\bX) \subset h^{\le 1}_{!*}(\bX) \subset h^{\le 2}_{!*}(\bX) 
\subset h^{\le 3}_{!*}(\bX) \subset h^{\le 4}_{!*}(\bX) = h_{!*}(\bX) \; .
\]
It is uniquely defined by the following property: both
the canonical projection
from $h(\Xp)$ to $h_{!*} (\bX)$ 
and the canonical inclusion of $h_{!*} (\bX)$ into $h(\Xp)$ 
are morphisms of filtered motives. 
The filtration is split in the sense that all 
$h^{\le r}_{!*}(\bX)$ admit direct
complements in $h_{!*}(\bX)$. In particular, the quotients
\[
h^{\ge r}_{!*}(\bX) := h_{!*}(\bX) / h^{\le {r-1}}_{!*}(\bX) 
\] 
of $h_{!*}(\bX)$ exist. \\[0.1cm]
(ii)~The filtration of $h_{!*}(\bX)$ is independent of the choice of
desingularization $\Xp$. \\[0.1cm]
(iii)~Duality $h_{!*}(\bX)^{\vee} \isoto h_{!*}(\bX)(2)$ 
(Proposition~\ref{2F}) induces isomorphisms
\[
h^{\le r}_{!*}(\bX)^{\vee} \isoto h^{\ge 4-r}_{!*}(\bX)(2) \; .
\]
\end{Thm}

\begin{Proof}
Define
\[
h^{\le r}_{!*}(\bX) := h^{\le r}(\Xp) \quad \text{for} \quad r \le 1
\]
and
\[
h^{\le r}_{!*}(\bX) := \ker (p^{\le r}) \quad \text{for} \quad r \ge 2 \; .
\]
Claim (i) is a consequence of the compatibility of the idempotents $p^{\le r}$,
(ii) is a consequence of Proposition~\ref{2E}~(iv),
and (iii) follows from symmetry of $p$.  
\end{Proof}

\begin{Def} \label{3e}
(a)~The filtration
\[
0 \subset h^0_{!*}(\bX) \subset h^{\le 1}_{!*}(\bX) 
\subset h^{\le 2}_{!*}(\bX) 
\subset h^{\le 3}_{!*}(\bX) \subset h^{\le 4}_{!*}(\bX) = h_{!*}(\bX) \; .
\]
from Theorem~\ref{3C} is called the \emph{K\"unneth filtration 
of $h_{!*}(\bX)$}. \\[0.1cm]
(b)~The $n$-th \emph{K\"unneth component of $h_{!*}(\bX)$}, $0 \le n \le 4$,
is the sub-quotient of $h_{!*}(\bX)$ defined by
\[
h^n_{!*}(\bX) := h^{\le n}_{!*}(\bX) / h^{\le n-1}_{!*}(\bX) \; .
\]
\end{Def}

For future reference, let us note the following immediate consequence of
our construction.

\begin{Prop} \label{3f}
Let $n$ be an integer unequal to two. Then there is a canonical isomorphism
of motives
\[ 
h^n_{!*}(\bX) \isoto h^n(\Xp) \; .
\]
\end{Prop}

\begin{Rem}
One may define the notion of K\"unneth decomposition of the intersection
motive as being a decomposition splitting the K\"unneth filtration.
Adding the complement $\oplus_m h^2 (D_m)$ of $h_{!*}(\bX)$ in $h(\Xp)$,
one gets a K\"unneth decomposition of $h(\Xp)$ in the abstract sense of
Remark~\ref{3c}~(b). 
It is not clear to me whether such a K\"unneth decomposition
does necessarily occur among those obtained using Murre's construction recalled earlier, when
$D$ has more than one component. The
problem is the relation
\[
p_{\Xp,*}(p_3 \cdot [\Xp \times D_m]) = 0 \; ;
\]
here as in the sequel, we use the same notation as in the proof of Proposition~\ref{3B}. 
The cycle class in question is
a non-zero multiple of
\[
p_{\Xp,*} (\wb \cdot [\Xp \times C \cdot D_m]) \; .
\]
For any fixed $m$, the K\"unneth decomposition of $h(\Xp)$ 
can be \emph{chosen} such that this cycle
class vanishes: take $T$ to be equal to $\frac{1}{d} [C \cdot D_m]$, where $d$ is the degree of $C
\cdot D_m$. 
\end{Rem}

%%% Local Variables:
%%% mode: latex
%%% TeX-master: "imot"
%%% End:

\bigskip

%\include{Sec4}
%%%%%%%%%%%%%%%%%%%%%%%%%%%%%%%%%%%%%%%%%%%%%%%%%%%%%%%%%%%%%%%%%%%%%%%
%
%  Section 4
%
%%%%%%%%%%%%%%%%%%%%%%%%%%%%%%%%%%%%%%%%%%%%%%%%%%%%%%%%%%%%%%%%%%%%%%%

\section{Hard Lefschetz for the intersection motive}
\label{4}

%%%%%%%%%%%%%%%%%%%%%%%%%%%%%%%%

%%%%%%%%%%%%%%%%%%%%%%%%%%%%%%%%

We continue to consider 
a proper surface $\bX$ over the
base field $k$. Let us consider the K\"unneth filtration
\[
0 \subset h^0_{!*}(\bX) \subset h^{\le 1}_{!*}(\bX) \subset 
h^{\le 2}_{!*}(\bX) \subset h^{\le
3}_{!*}(\bX) \subset h^{\le 4}_{!*}(\bX) = h(\bX)_{!*}
\]
of the intersection motive. The aim of this section is to prove the following.

\begin{Thm} \label{4A}
Let $\CL$ be a line bundle on $\bX$. \\[0.1cm]
(i)~There is a morphism of motives
\[
c_{\CL}: h_{!*}(\bX)(-1) \longto h_{!*}(\bX) \; ,
\]
which is uniquely characterized by the following two properties: 
\begin{enumerate}
\item[(1)] If $\bX$ is smooth, then $c_{\CL}$ equals the
cup-product with the first Chern class of $\CL$ on 
$h(\bX)(-1) = h_{!*}(\bX)(-1)$
\cite[Sect.~2.1]{Sch}. 
\item[(2)] The morphism
$c_{\CL}$ is contravariantly functorial with respect to dominant morphisms
$g: \bY \to \bX$ of proper surfaces over $k$: the diagram
\[
\vcenter{\xymatrix@R-10pt{
        h_{!*}(\bY)(-1) \ar[r]^-{c_{g^* \! \CL}} &
        h_{!*}(\bY) \\
        h_{!*}(\bX)(-1) \ar[r]^-{c_{\CL}} \ar[u]^{g^*(-1)} &
        h_{!*}(\bX) \ar[u]_{g^*}
\\}}
\]
(see Proposition~\ref{2E}~(i)) commutes. 
\end{enumerate}
(ii)~If $\CL'$ is a second line bundle on $\bX$, then 
\[
c_{\CL \otimes \CL'} = c_{\CL} + c_{\CL'} \; .
\]
In other words, the map
\[
\Pic(\bX) \longto \Hom \bigl( h_{!*}(\bX)(-1), h_{!*}(\bX) \bigr) \; , \;
\CL \longmapsto c_{\CL}
\]
is a morphism of groups. \\[0.1cm]
(iii)~The morphism $c_{\CL}$ is filtered in the following sense: it
induces morphisms
\[
c_{\CL}: h^{\le n-2}_{!*}(\bX)(-1) \longto h^{\le n}_{!*}(\bX) 
\]
and hence, morphisms
\[
c_{\CL}: h^{n-2}_{!*}(\bX)(-1) \longto h^n_{!*}(\bX) 
\]
for all $n \in \BZ$. \\[0.1cm]
(iv)~If ($\bX$ is projective and) $\CL$ or $\CL^{-1}$ is ample, then
\[
c_{\CL}^2 = c_{\CL} \circ c_{\CL}: h^0_{!*}(\bX)(-2) \longto h^4_{!*}(\bX) 
\] 
and
\[
c_{\CL}: h^1_{!*}(\bX)(-1) \longto h^3_{!*}(\bX) 
\]
are isomorphisms.
\end{Thm}

Part (iv) of this result should be seen as the motivic analogue of
the Hard Lefschetz Theorem for intersection cohomology 
\cite[Thm.~6.2.10]{BBD}. \\

In order to prepare the proof of Theorem~\ref{4A}, let us recall the
ingredients of the proof when $\bX$ is smooth
(in which case Theorem~\ref{4A} is of course known).  
The morphism $c_{\CL}$ then equals the cup-product
with the first Chern class, which can be described as follows. 
In the category $CHM(k)_{\BQ}$, the vector space $CH^1(\bX)$
equals the group of morphisms from $\BL$ to $h(\bX)$.
We define $c_{\CL}$ as being the composition
\[
h(\bX)(-1) = h(\bX) \otimes \BL \stackrel{\id_{\bX}^* \otimes [\CL]}{\longto}
h(\bX) \otimes h(\bX) \stackrel{\Delta^*}{\longto} h(\bX)
\]
($\Delta:= $ the diagonal embedding $\bX \into \bX \times_k \bX$).
From this description, pro\-perties (i)~(2) (for smooth $\bY$) and (ii)
are immediate. Recall that $\bX$, as a smooth and proper
surface, is projective.
Since the group $\Pic(\bX)$ is generated by the classes
of very ample line bundles, in order to prove (iii) and (iv), we may
(by (ii)) assume that $\CL$ is very ample. In addition, we may prove the
claims after base change to a finite extension of $k$, and hence assume
that $\bX$ is geometrically connected, and that $\CL$ is represented by 
a smooth curve $C$
embedded into $\bX$ via the closed immersion $i_C$. The morphism
$c_{\CL}$ then equals the composition of
\[
i_C^*(-1): h(\bX)(-1) \longto h(C)(-1)
\]
and of
\[
i_{C,*}: h(C)(-1) \longto h(\bX) \; .
\]
By auto-duality of the K\"unneth filtrations for $C$ and for $\bX$,
it suffices for (iii) to show that $i_C^*: h(\bX) \to h(C)$ 
is a morphism of filtered
motives. But this follows from \cite[Lemma~8.3.2]{Mr3} and 
\cite[Prop.~5.8]{J}.
As for (iv), observe that identifying $h^0(\Xp)(-2)$ and $h^4(\Xp)$
with $\BQ(-2)$ allows to relate the morphism 
$c_{\CL}^2: h^0(\Xp)(-2) \to h^4(\Xp)$
to the self-intersection number $C \cdot C$,
which is strictly positive since $\CL$ is very ample.
The statement on
$c_{\CL}: h^1(\Xp)(-1) \to h^3(\Xp)$ is the most difficult to prove.
We refer to \cite[Thm.~4.4~(ii)]{Sch} for the details. \\

Given the contravariance property of the intersection motive
(Proposition~\ref{2E}~(i)), it is now clear what remains to be done
in order to prove Theorem~\ref{4A} in the generality we stated it. 
First note that in our statement, we may replace $\bX$ by its
normalization $\Xs$.
Indeed, $h_{!*}(\bX) = h_{!*}(\Xs)$, and the morphism
$\Xs \to \bX$ being finite, the pull-back of an ample line bundle on $\bX$
is ample on $\Xs$. 
Next, fix a cartesian diagram
\[
\vcenter{\xymatrix@R-10pt{
        X \ar@{^{ (}->}[r] \ar@{=}[d] &
        \Xp \ar@{<-^{ )}}[r]^{\ip} \ar[d]_\pi &
        D \ar[d]^\pi \\
        X \ar@{^{ (}->}[r] &
        \Xs \ar@{<-^{ )}}[r] &
        Z
\\}}
\]
which is a desingularization of $\Xs$.
Thus, $\pi$ is proper, $\Xp$ is smooth and proper (hence projective), 
$Z$ is finite,
and $D$ a divisor with normal crossings,
whose irreducible components $D_m$ are smooth. We need to show
that for any line bundle $\CL$ on $\Xs$,
the composition
\[
h_{!*}(\bX)(-1) \longinto h(\Xp)(-1) \stackrel{c_{\pi^* \! \CL}}{\longto} h(\Xp)
\]
lands in $h_{!*}(\bX) \subset h(\Xp)$ --- this will then be our definition
of $c_{\CL}$ --- and that we have the Hard Lefschetz Theorem \ref{4A}~(iv). 
In fact, we shall prove a more general result.

\begin{Var} \label{4A'}
Let 
$\widetilde{\CL}$ be a line bundle 
on $\Xp$, whose restrictions to all $D_m$ are trivial 
(for example, the pull-back of a line bundle on $\Xs$). \\[0.1cm]
(i)~The restriction of the morphism of motives
\[
c_{\widetilde{\CL}} : h(\Xp)(-1) \longto h(\Xp) 
\] 
to the sub-motive $h_{!*}(\bX)(-1)$  
induces a morphism $h_{!*}(\bX)(-1) \to h_{!*}(\bX)$. 
In other words,
there is a commutative diagram
\[
\vcenter{\xymatrix@R-10pt{
        h(\Xp)(-1) \ar[r]^-{c_{\widetilde{\CL}}} &
        h(\Xp) \\
        h_{!*}(\bX)(-1) \ar[r]^-{c_{\widetilde{\CL}}} 
                                        \ar@{_{ (}->}[u]^{\pi^*(-1)} &
        h_{!*}(\bX) \ar@{^{ (}->}[u]_{\pi^*}
\\}}
\]
(ii)~If $\widetilde{\CL}'$ is a second line bundle on $\Xp$
with trivial restrictions to all $D_m$, 
then 
\[
c_{\widetilde{\CL} \otimes \widetilde{\CL}'} 
= c_{\widetilde{\CL}} + c_{\widetilde{\CL}'} \; .
\]
(iii)~The morphism $c_{\widetilde{\CL}}$ is filtered: it
induces morphisms
\[
c_{\widetilde{\CL}}: h^{\le n-2}_{!*}(\bX)(-1) \longto h^{\le n}_{!*}(\bX) 
\]
for all $n \in \BZ$. \\[0.1cm]
(iv)~Assume in addition
that $\widetilde{\CL}$ is the line bundle
associated to a divisor 
$C$ on $\Xp$ 
such that
$C - \sum_m a_m D_m$ or $-C - \sum_m a_m D_m$
is ample for a suitable choice of integers $a_m \ge 0$
(for example, $\widetilde{\CL} = \pi^* \! \CL$ for an ample line
bundle  $\CL$ on $\Xs$). Then
\[
c_{\widetilde{\CL}}^2: h^0_{!*}(\bX)(-2) \longto h^4_{!*}(\bX) 
\] 
and
\[
c_{\widetilde{\CL}}: h^1_{!*}(\bX)(-1) \longto h^3_{!*}(\bX) 
\]
are isomorphisms.
\end{Var}

\begin{Proof}
In order to prove (i), we have to check that the composition
\[
h_{!*}(\bX)(-1) \stackrel{\pi^*(-1)}{\longinto} h(\Xp)(-1) 
\stackrel{c_{\widetilde{\CL}}}{\longto} h(\Xp) 
\stackrel{{\ip}_* \alpha^{-1} {{\ip}}^*}{\longto} h(\Xp)
\]
is zero. Since the formation of Chern classes is compatible with
pull-backs, the composition ${{\ip}}^* c_{\widetilde{\CL}}$ equals
\[
h(\Xp)(-1) \stackrel{\oplus_m i_m^*}{\longto} \bigoplus_m h(D_m)(-1)
\stackrel{\oplus_m c_{i_m^* \! \widetilde{\CL}}}{\longto}
\bigoplus_m h(D_m) \longonto \bigoplus_m h^2(D_m) \; ,
\] 
where $i_m$ denotes the immersion of $D_m$ into $\Xp$. But by assumption,
the morphisms
$c_{i_m^* \! \widetilde{\CL}}: h(D_m)(-1) \to h(D_m)$ are all zero.
 
Claims (ii) and (iii) hold since they
hold for $c_{\widetilde{\CL}} : h(\Xp)(-1) \to h(\Xp)$.

As for (iv), observe that according to Proposition~\ref{3f}, 
\[ 
h^n_{!*}(\bX) \cong h^n(\Xp) \; , \; n \ne 2 \; .
\]
Thus, we have to prove that
\[
c_{\widetilde{\CL}}^2: h^0(\Xp)(-2) \longto h^4(\Xp) 
\] 
and
\[
c_{\widetilde{\CL}}: h^1(\Xp)(-1) \longto h^3(\Xp) 
\]
are isomorphisms. As before, the claim for $c_{\widetilde{\CL}}^2$
is essentially equivalent to showing that the self-intersection number
$C \cdot C$ is non-zero. Since the restriction of $\widetilde{\CL}$
to any of the $D_m$ is trivial, 
we have the formula
\[
C \cdot C = 
\bigl( \pm C - \sum_m a_m D_m \bigr) \cdot 
\bigl( \pm C - \sum_m a_m D_m \bigr) -
\bigl( \sum_m a_m D_m \bigr) \cdot \bigl( \sum_m a_m D_m \bigr) \; .
\]
The intersection matrix $(D_n \cdot D_m)_{n,m}$
is negative definite \cite[p.~6]{M}, hence the matrix
$\bigl( (a_n D_n) \cdot (a_m D_m) \bigr)_{n,m}$ is negative semi-definite.
It follows that the term 
$(\sum_m a_m D_m) \cdot (\sum_m a_m D_m)$ is non-positive. Hence
\[
C \cdot C \ge 
\bigl( \pm C - \sum_m a_m D_m \bigr) \cdot 
\bigl( \pm C - \sum_m a_m D_m \bigr) \; .
\]
But by assumption,
one of the divisors $C - \sum_m a_m D_m$, $-C - \sum_m a_m D_m$ is ample.
Therefore, its self-intersection number 
is strictly positive. 

In order to prove the claim for 
$c_{\widetilde{\CL}}: h^1(\Xp)(-1) \to h^3(\Xp)$, 
observe first that by (ii), we may assume $C - \sum_m a_m D_m$ to be very ample.
By passing to a finite extension of $k$, we find a smooth curve $H$ embedded
into $\Xp$ via the closed immersion $i_H$, and such that 
there is an equivalence of divisors
\[
C - \sum_m a_m D_m \sim H \; .
\] 
In particular, $H$ is very ample, and
\[
c_{\widetilde{\CL}} =  i_{H,*} i_H^* + \sum_m a_m i_{m,*} i_m^* :
h^1(\Xp)(-1) \longto h^3(\Xp) \; .
\]
Hard Lefschetz \ref{4A}~(iv) tells us that $i_{H,*} i_H^*$ is
an isomorphism. In order to see that the same still holds after
adding the ``error term'' $\sum_m a_m i_{m,*} i_m^*$, we neeed to  
recall more details of the proof.

In fact, as follows from \cite[Prop.~4.5]{Sch}, the full sub-category
of motives isomorphic to $h^1(Y)$, for smooth projective varieties
$Y$ over $k$, is equivalent to the category of Abelian varieties over $k$
up to isogeny. More precisely, this equivalence is such that
$h^1(Y)$ corresponds to the Picard variety $P_Y$, and that the motive
$h^{2d_Y-1}(d_Y-1)$ (for $Y$ of pure dimension $d_Y$) corresponds to the
Albanese variety $A_Y$. Furthermore, for a morphism $f:Y_1 \to Y_2$,
the pull-back of motives $f^*: h^1(Y_2) \to h^1(Y_2)$ corresponds to
$f^*: P_{Y_2} \to P_{Y_1}$, while the push-forward
$f_*: h^{2d_{Y_1}-1}(d_{Y_1}-1) \to h^{2d_{Y_2}-1}(d_{Y_2}-1)$
(for $Y_i$ of pure dimension $d_{Y_i}$, $i=1,2$) corresponds to
$f_*: A_{Y_1} \to A_{Y_2}$. Proving that $c_{\widetilde{\CL}}$ is
an isomorphism of motives is thus equivalent to proving the following statement:
the composition of
\[
I^* : P_{\widetilde{X}} \longto P_H \times_k \prod_m \bigl( P_{D_m} \bigr)^{a_m}
\] 
with its dual
\[
I_* : A_H \times_k \prod_m \bigl( A_{D_m} \bigr)^{a_m} \longto A_{\widetilde{X}}
\]
is an isogeny from the Picard variety of $\widetilde{X}$ to the
Albanese variety of $\widetilde{X}$ (recall that our motives are with $\BQ$-coefficients). 
Here, $I$ denotes the morphism from 
the disjoint union of $H$ and $a_m$ copies of $D_m$, for all $m$, to
$\widetilde{X}$. Also, we have identified the Picard and the Albanese varieties
of the curves $H$ and $D_m$ to the respective Jacobians,
using the fact that these are canonically principally polarized.

The decisive ingredient of the proof is \cite[Cor.~1 of Thm.~7]{We}, which states
that since $H$ is very ample, the kernel of
$i_H^*: P_{\widetilde{X}} \to P_H$ is finite. The same is thus true
for $I^*$. Now observe that a polarization on an Abelian variety 
(such as $P_H \times_k \prod_m \bigl( P_{D_m} \bigr)^{a_m}$) induces a
polarization on any sub-Abelian variety. The composition
$I_* I^*$ is therefore an isogeny.
\end{Proof}

%%% Local Variables:
%%% mode: latex
%%% TeX-master: "imot"
%%% End:

\bigskip

%\include{Sec5}
%%%%%%%%%%%%%%%%%%%%%%%%%%%%%%%%%%%%%%%%%%%%%%%%%%%%%%%%%%%%%%%%%%%%%%%
%
%  Section 5
%
%%%%%%%%%%%%%%%%%%%%%%%%%%%%%%%%%%%%%%%%%%%%%%%%%%%%%%%%%%%%%%%%%%%%%%%

\section{The motive of the exceptional divisor}
\label{5}

%%%%%%%%%%%%%%%%%%%%%%%%%%%%%%%%

%%%%%%%%%%%%%%%%%%%%%%%%%%%%%%%%

At this point, we need to enlarge the category of motives we are working in
since we wish to 
consider motives of genuinely singular varieties. Let us first set up
the notation, which 
follows that of \cite{VSF}. From now on, our base field $k$
is assumed to be perfect. 
We write $Sch/k$ for the category of schemes which are
separated and of finite type
over $k$, and $Sm/k$ for the full sub-category of objects of $Sch/k$ which
are smooth over $k$. Recall the definition
of the category $\SmC$ \cite[p.~190]{VSF}:
its objects are those of $Sm/k$. Morphisms
from $Y$ to $X$ are given by the group $c(Y,X)$ of \emph{finite
correspondences} from $Y$ to $X$.  The category $\ShN$
of \emph{Nisnevich sheaves with transfers}
\cite[Def.~3.1.1]{VSF} is the category of those
contravariant additive functors from $\SmC$ to Abelian groups,
whose restriction to $Sm/k$ is a sheaf for the Nisnevich topo\-logy.
Inside the derived category $D^-(\ShN)$ of complexes bounded from
above, one defines the full triangulated sub-category $\DM$
of \emph{effective motivic complexes} over $k$
\cite[p.~205, Prop.~3.1.13]{VSF} as the one consisting
of objects whose cohomology sheaves are \emph{homotopy invariant}
\cite[Def.~3.1.10]{VSF}.
The inclusion of
$\DM$ into $D^-(\ShN)$ admits a left adjoint $\RC$, which is
induced from the functor
\[
\uC_*: \ShN \longto C^-(\ShN) \; .
\]
By definition, $\uC_*$ maps $F \in \ShN$
to the \emph{singular simplicial complex} $\uC_*(F)$ \cite[p.~207, Prop.~3.2.3]{VSF}. 
One defines a functor $L$ from $Sch/k$ to
$\ShN$: it associates
to $X$ the Nisnevich sheaf with transfers $c(\argdot,X)$;
note that the above definition of $c(Y,X)$ still makes
sense when $X \in Sch/k$ is not
necessarily smooth. One defines
the \emph{motive} $M (X)$ as $\RC (L(X))$.
We shall use the same symbol for $M (X) \in \DM$ and for
its canonical representative $\uC_* (L(X))$ in $C^- (\ShN)$.
There is a second functor $L^c$, which associates to $X \in Sch/k$
the Nisnevich sheaf of quasi-finite correspondences
\cite[p.~223, 224]{VSF}.
One defines the \emph{motive with compact support}
$M^c (X)$ of $X \in Sch/k$ as
$\RC (L^c(X))$. 
It coincides with $M(X)$ if $X$ is proper. \\

A second, more geometric approach to motives is
the one developed in \cite[Sect.~2.1]{VSF}. 
There, the triangulated
category $\DeffgM$ of \emph{effective geometrical motives} over $k$
is defined. There is a canonical
full triangulated embedding of $\DeffgM$ into $\DM$ \cite[Thm.~3.2.6]{VSF},
which maps the geometrical
motive of $X \in Sm/k$ \cite[Def.~2.1.1]{VSF} to $M (X)$.
Using this embedding, we consider $M (X)$ as an object of $\DeffgM$.
The \emph{Tate motive} $\BZ(1)$ in $\DeffgM$ is defined as the \emph{reduced
motive} of $\BP^1_k$ \cite[p.~192]{VSF}, shifted by $-2$. 
There is a canonical direct sum decomposition
\[
M(\BP^1_k) = \BZ(0) \oplus \BZ(1)[2] \; .
\]
The category $\DgM$ of \emph{geometrical motives} over $k$
is obtained from the category $\DeffgM$ by inverting $\BZ(1)$.
All categories $\DeffgM$, $\DgM$, $D^-(\ShN)$, and $\DM$
are tensor triangulated, and admit unit objects, which we denote
by the same symbol $\BZ (0)$
\cite[Prop.~2.1.3, Cor.~2.1.5, p.~206,
Thm.~3.2.6]{VSF}.
For $M \in \DgM$ and $n \in \BZ$, write $M(n)$ for the tensor product
$M \otimes \BZ (n)$.
According to \cite[Cor.~4.10]{V3},
the functor $\DeffgM \to \DgM$ is a full triangulated
embedding (see \cite[Thm.~4.3.1]{VSF} for a proof when $k$ admits
resolution of singularities). \\

Let us denote by $\DeffgM_{\BQ}$ and $\DgM_{\BQ}$ the triangulated categories
obtained by the
$\BQ$-linear analogues of the above constructions \cite[Sect.~16.2.4
and Sect.~17.1.3]{A}. The relation
to Chow motives is given by the following result due to Voevodsky.

\begin{Thm} \label{5A}
(i)~There is a natural contravariant $\BQ$-linear tensor functor
\[
R : CHM(k)_{\BQ} \longto \DgM_{\BQ} \; .
\]
$R$ is fully faithful. \\[0.1cm]
(ii)~For any smooth projective variety $S$ over $k$,
the functor $R$ maps the Chow motive $h(S)$ to the motive 
$M(S) \in \DeffgM_{\BQ} \subset \DgM_{\BQ}$. \\[0.1cm]
(iii)~The functor $R$ maps the Lefschetz motive $\BL$ to the motive 
$\BZ(1)[2]$, compatibly with the decompositions
\[
h(\BP^1_k) = h(\Spec k) \oplus \BL
\]
in $CHM(k)_{\BQ}$ and
\[
M(\BP^1_k) = \BZ(0) \oplus \BZ(1)[2] 
\]
in $\DeffgM_\BQ$. 
\forget{
(iv)~If $k$ admits resolution of singularities,
thenthere are no non-trivial higher extensions
between Chow motives in the category $\DgM_{\BQ}$. 
More precisely, given smooth projective varieties
$S_1$, $S_2$ over $k$ and $i > 0$, the group
\[
\Hom_{\DgM_{\BQ}} \bigl( M(S_1) , M(S_2)[i] \bigr)
\]
is zero.
}
\end{Thm}

\begin{Proof}
The essential point of the proof is to show equality of morphisms:
\[
\Hom_{CHM(k)_{\BQ}} \bigl( h(Y)(-q) , h(X) \bigr) = 
\Hom_{\DgM_{\BQ}} \bigl( M(X) , M(Y)(q)[2q] \bigr) 
\]
for smooth projective varieties $X$ and $Y$ over $k$ and $q \ge 0$. 
Duality in $\DgM_{\BQ}$
\cite[Thm.~18.4.1.1]{A} (\cite[Thm.~4.3.7]{VSF} if $k$ admits
resolution of singularities) allows us to reduce to the case
$Y = \Spec k$, in which case the claim follows from \cite[Cor.~2]{V2}.
\end{Proof}

\forget{
\begin{Rem}
It should be noted that the preceding result holds already before tensoring
with $\BQ$, i.e., there is a version for Chow motives with $\BZ$-coefficients.
But we intend to apply the results of the preceding sections, where 
invertibility of non-zero integers was needed. First, the determinant $d$
of the intersection matrix $(D_n \cdot D_m)_{n,m}$ is a non-zero integer, a fact
which is essential for the construction of the intersection motive
(Theorem~\ref{2B}~(i)). This could be controlled by tensoring only with
$\BZ[1/d]$. Second, and more seriously, our proofs systematically made use
of the principle that a finite extension of the base field does not affect
the statement in question.
This principle is of course based on the invertibility of the degree
of the field extension. In general, it will not be possible to 
control these degrees (unless for example, our field is algebraically closed).  
\end{Rem}
}

\begin{Ex} \label{5B}
Fix a proper surface $\bX$ over $k$. Recall the K\"unneth filtration
of the intersection motive
\[
0 \subset h^0_{!*}(\bX) \subset h^{\le 1}_{!*}(\bX) 
\subset h^{\le 2}_{!*}(\bX) 
\subset h^{\le 3}_{!*}(\bX) \subset h^{\le 4}_{!*}(\bX) = h_{!*}(\bX) \; ,
\]
the quotients
\[
h^{\ge r}_{!*}(\bX) := h_{!*}(\bX) / h^{\le r-1}_{!*}(\bX) \; ,
\]
and the K\"unneth components
\[
h^n_{!*}(\bX) = h^{\le n}_{!*}(\bX) / h^{\le n-1}_{!*}(\bX) 
\]
(Definition~\ref{3e}). Let us write 
$M^{!*}(\bX) := R(h_{!*}(\bX))$, 
\[
M^{!*}_{\ge r}(\bX) := R(h^{\ge r}_{!*}(\bX)) \; ,
\]
\[
M^{!*}_{\le n}(\bX) := R(h^{\le n}_{!*}(\bX)) \; ,
\]
\[
M^{!*}_n(\bX) := R(h^n_{!*}(\bX)) \; .
\]
We thus have exact triangles
\[
M^{!*}_{\ge r+1}(\bX) \longto M^{!*}(\bX) \longto 
M^{!*}_{\le r}(\bX) \stackrel{\delta}{\longto} M^{!*}_{\ge r+1}(\bX)[1] \; ,
\] 
\[
M^{!*}_n(\bX) \longto M^{!*}_{\le n}(\bX) \longto 
M^{!*}_{\le n-1}(\bX) \stackrel{\delta}{\longto} M^{!*}_n(\bX)[1] 
\] 
in $\DeffgM_{\BQ}$, which are all split in the sense that the 
boundaries $\delta$
are zero.
\end{Ex}

For the rest of this section, fix a (not necessarily proper)
surface $\bX$ over $k$, and a cartesian diagram
\[
\vcenter{\xymatrix@R-10pt{
        X \ar@{^{ (}->}[r] \ar@{=}[d] &
        \Xp \ar@{<-^{ )}}[r] \ar[d]_\pi &
        D \ar[d]^\pi \\
        X \ar@{^{ (}->}[r] &
        \Xs \ar@{<-^{ )}}[r] &
        Z
\\}}
\]
which is a desingularization of the normalization $\Xs$.
Thus, $\pi$ is proper, $\Xp$ is smooth, $Z$ is finite,
and $D$ a divisor with normal crossings,
whose irreducible components $D_m$ are smooth projective curves. 
The exact triangle associated to the closed covering of
$D$ by the $D_m$ \cite[Prop.~4.1.3]{VSF} (but see also the proof 
of Proposition~\ref{6D}~(i))
shows that $M(D)$ belongs 
to the category $\DeffgM$.

\begin{Def} \label{5C}
Define Chow motives $h^0(D)$ and $h^2(D)$ as follows. \\[0.1cm]
(a)~$h^0(D) := h(S)$, where $S$ equals the spectrum of the ring of
global sections of the structure sheaf of $D$. \\[0.1cm]
(b)~$h^2(D) := \oplus_m h^2 (D_m)$.
\end{Def}

Let us write $M_0(D) := R(h^0(D))$ and $M_2(D) := R(h^2(D))$.
The morphism $D \to S$ and the inclusions $i_m$ of the components
$D_m$ into $D$ induce morphisms $M(D) \to M_0(D)$
and $M_2(D) \to M(D)$ in $\DeffgM_{\BQ}$.

\begin{Lem} \label{5D}
The morphism $M(D) \to M_0(D)$ is a split epimorphism, and 
$M_2(D) \to M(D)$ is a split monomorphism.
The composition of the two morphisms $M_2(D) \to M(D) \to M_0(D)$ is trivial.
\end{Lem}

\begin{Proof}
The composition
\[
\bigoplus_m R(h^0(D_m)) \longto \bigoplus_m R(h(D_m)) = \bigoplus_m M(D_m)
\longto M(D) \longto M_0(D)
\]
is a split epimorphism, hence so is $M(D) \to M_0(D)$. The composition
\[
M_2(D) \longto M(D) \longto M(\Xp)
\]
is a split monomorphism (Theorem~\ref{2B}~(i)),
hence so is $M_2(D) \to M(D)$.
The last claim is obvious.
\end{Proof}

It follows that the objects 
\[
M_{\ge 1}(D) := \ker \bigl( M(D) \longto M_0(D) \bigr) \; ,
\]
\[
M_{\le 1}(D) := M(D) / M_2(D) \; ,
\]
and  
\[
M_1(D) := \ker \bigl( M_{\le 1}(D) \longto M_0(D) \bigr)
= M_{\ge 1}(D) / M_2(D)
\]
exist.
They give rise to what we might call the K\"unneth filtration of 
$M(D)$:
\[
M(D) =: M_{\le 2}(D) \longonto M_{\le 1}(D) \longonto M_0(D) \; ,
\]
\[
M_2(D) \longinto M_{\ge 1}(D) \longinto M_{\ge 0}(D) := M(D) \; .
\]
Note that there are split exact triangles
\[
M_2(D) \longto M(D) \longto 
M_{\le 1}(D) \stackrel{\delta = 0}{\longto} M_2(D)[1] \; ,
\]
\[
M_1(D) \longto M_{\le 1}(D) \longto 
M_0(D) \stackrel{\delta = 0}{\longto} M_1(D)[1] 
\] 
in $\DeffgM_{\BQ}$.
For all $m$,
let us also define $M_i(D_m)$, $0 \le i \le 2$ and $M_{\le 1} (D_m)$
as the images under the functor $R$ of the Chow motives
$h^i(D_m)$ and $h^{\le 1} (D_m)$, respectively.

\begin{Rem}
Unlike $M_0(D)$ and $M_2(D)$, the sub-quotient $M_1(D)$ should
not in general be expected to come from a Chow motive. Indeed, as
we shall see, the ``kernel'' of 
\[
\bigoplus_{n < m} M(D_n \cap D_m)[1] \longto \bigoplus_m M_0(D_m)[1]
\]
contributes to $M_1(D)$.  
\end{Rem}

\bigskip

%\include{Sec6}
%%%%%%%%%%%%%%%%%%%%%%%%%%%%%%%%%%%%%%%%%%%%%%%%%%%%%%%%%%%%%%%%%%%%%%%
%
%  Section 6
%
%%%%%%%%%%%%%%%%%%%%%%%%%%%%%%%%%%%%%%%%%%%%%%%%%%%%%%%%%%%%%%%%%%%%%%%

\section{An extension of motives}
\label{6}

%%%%%%%%%%%%%%%%%%%%%%%%%%%%%%%%

%%%%%%%%%%%%%%%%%%%%%%%%%%%%%%%%

We continue to study the situation
\[
\vcenter{\xymatrix@R-10pt{
        X \ar@{^{ (}->}[r]^{\jp} \ar@{=}[d] &
        \Xp \ar@{<-^{ )}}[r]^{\ip} \ar[d]_\pi &
        D \ar[d]^\pi \\
        X \ar@{^{ (}->}[r]^j &
        \Xs \ar@{<-^{ )}}[r] &
        Z
\\}}
\]
fixed in Section~\ref{5}, but assume in addition that the surface $\bX$
is proper. The morphism $\ip_*: M(D) \to M(\Xp)$ will be at the
base of the construction of an extension in 
$\DeffgM_{\BQ}$ (Theorem~\ref{Main}).
Let us start with a number of elementary observations.

\begin{Lem} \label{6a}
The composition 
\[
M(D) \stackrel{\ip_*}{\longto} M(\Xp) \longonto M^{!*}(\bX)
\]
factors uniquely through a morphism $\ip_*: M_{\le 1}(D) \to M^{!*}(\bX)$.
\end{Lem}

\begin{Proof}
We identify $M^{!*}(\bX)$ with the categorical quotient of $M(\Xp)$ by
$M_2(D)$. The composition in question
thus vanishes on $M_2(D)$. It therefore factors uniquely over 
the categorical quotient $M_{\le 1}(D)$ 
of $M(D)$ by $M_2(D)$. 
\end{Proof}

\begin{Rem} \label{6B}
If $k$ admits resolution of singularities,
then we have \emph{localization} for the motive with compact support 
\cite[Prop.~4.1.5]{VSF}. In our situation, this means that there is a
canonical exact triangle
\[
M(D) \stackrel{\ip_*}{\longto} M(\Xp) 
\stackrel{\jp^*}{\longto} M^c(X) \longto M(D)[1] \; .
\]
From this, one deduces easily that 
$\ip_*: M_{\le 1}(D) \to M^{!*}(\bX)$ sits in an exact triangle
\[
M_{\le 1}(D) \stackrel{\ip_*}{\longto} M^{!*}(\bX) 
\stackrel{j^*}{\longto} M^c(X) \longto M_{\le 1}(D)[1] \; .
\]
\end{Rem}

Consider the sub-object $M_1(D)$ of $M_{\le 1}(D)$, 
and the quotient $M_0^{!*}(\bX)$ of $M^{!*}(\bX)$.

\begin{Lem} \label{6C}
The composition
\[
M_1(D) \longinto M_{\le 1}(D) \stackrel{\ip_*}{\longto} 
M^{!*}(\bX) \longonto M_0^{!*}(\bX)
\]
is trivial.
\end{Lem}

\begin{Proof}
The motive $M_0^{!*}(\bX)$ equals $M_0(\Xp) := R(h^0(\Xp))$
(Proposition~\ref{3f}), hence the composition
\[
M_{\le 1}(D) \stackrel{\ip_*}{\longto} 
M^{!*}(\bX) \longonto M_0^{!*}(\bX)
\]
equals the composition
\[
M_{\le 1}(D) \longonto M_0(D) \stackrel{\ip_*}{\longto} M_0(\Xp) \; .
\]
It is therefore trivial on $M_1(D)$.
\end{Proof}

\begin{Cor} \label{6Ca}
The morphism $\ip_*: M_{\le 1}(D) \to M^{!*}(\bX)$
respects the K\"un\-neth filtrations. 
\end{Cor}

The inclusion $\ip$ therefore induces
a morphism, equally denoted 
$\ip_*$ from $M_1(D)$ to $M^{!*}_{\ge 1}(\bX)$.
Consider the quotient $M_1^{!*}(\bX)$ of $M^{!*}_{\ge 1}(\bX)$.

\begin{Prop} \label{6D}
Assume that all geometric 
irreducible components of $D$ are of genus zero.
\\[0.1cm]
(i)~The object $M_1(D)[-1]$ of $\DeffgM_{\BQ}$
is an Artin motive, i.e., it is isomorphic
to the motive of some zero-dimensional variety over $k$. 
More precisely, there is a canonical exact sequence
of Artin motives 
\[
0 \longto M_1(D)[-1] \longto \bigoplus_{n < m} M(D_n \cap D_m) 
\longto \bigoplus_m M_0(D_m) \; ,
\]
and $M_1(D)[-1]$ is a direct summand of $\oplus_{n < m} M(D_n \cap D_m)$. 
\\[0.1cm]
(ii)~The composition
\[
M_1(D) \stackrel{\ip_*}{\longto} 
M^{!*}_{\ge 1}(\bX) \longonto M_1^{!*}(\bX)
\]
is trivial.
\end{Prop}

\begin{Proof}
(i)~Consider the closed covering of $D$ by the $D_m$. It induces an
exact sequence of Nisnevich sheaves with transfers
\[
0 \longto \bigoplus_{n < m} L(D_n \cap D_m) \longto \bigoplus_m L(D_m) 
\longto L(D) \longto 0 \; , 
\]
hence an exact triangle
\[
\bigoplus_{n < m} M(D_n \cap D_m) \longto \bigoplus_m M(D_m) 
\longto M(D) \longto \bigoplus_{n < m} M(D_n \cap D_m)[1] \; .
\]
Given the definition of $M_2$, we get an exact triangle
\[
\bigoplus_{n < m} M(D_n \cap D_m) \longto \bigoplus_m M_{\le 1}(D_m) 
\longto M_{\le 1}(D) \longto \bigoplus_{n < m} M(D_n \cap D_m)[1] \; .
\]
But the $M_1(D_m)$ are zero by assumption. 
Hence the exact triangle takes the form
\[
\bigoplus_{n < m} M(D_n \cap D_m) \longto \bigoplus_m M_0(D_m) 
\longto M_{\le 1}(D) \longto \bigoplus_{n < m} M(D_n \cap D_m)[1] \; ;
\]
it thus belongs to the full triangulated sub-category $d_{\le 0} \DeffgM_{\BQ}$
gene\-rated by motives of dimension $0$. This triangulated
sub-category is canonically equivalent to the bounded derived category 
of the Abelian semi-simple category of Artin motives (with 
$\BQ$-coefficients) over $k$
\cite[Prop.~3.4.1 and Remark~2 following it]{VSF}. The sequence
\[
\bigoplus_{n < m} M(D_n \cap D_m) \longto \bigoplus_m M_0(D_m) 
\longto M_0(D) \longto 0
\]
of Artin motives is exact. From this and the above exact triangle, we see that
$M_1(D)[-1]$ is an Artin motive, which fits into an exact sequence
\[
0 \longto M_1(D)[-1] \longto \bigoplus_{n < m} M(D_n \cap D_m) 
\longto \bigoplus_m M_0(D_m) \; .
\]
(ii)~The motive $M_1^{!*}(\bX)$ equals $M_1(\Xp)$
(Proposition~\ref{3f}).
We shall show triviality of
\[
\Hom_{\DeffgM_{\BQ}} \bigl( M(Y)[1] , M_1(\Xp) \bigr)  
\]
for any smooth variety $Y$ over $k$. 
Applied to $Y = D_n \cap D_m$, $n < m$, together with (i), this will
establish (ii). Hard Lefschetz 
\[
M_1(\Xp) \cong M_3(\Xp)(-1)[-2]
\]
and duality
in $\DgM_{\BQ}$ imply that $\Hom_{\DeffgM_{\BQ}} ( M(Y)[1] , M_1(\Xp) )$ is isomorphic to
\[
\Hom_{\DgM_{\BQ}} \bigl( M_1(\Xp) \otimes M(Y)(-1)[-1] , \BZ(0) \bigr) \; ,
\]
which equals the direct factor 
$\Hom_{\DeffgM_{\BQ}} \bigl( M_1(\Xp) \otimes M(Y) , \BZ(1)[1] \bigr)$ of
\[
\Hom_{\DeffgM_{\BQ}} \bigl( M(\Xp \times_k Y) , \BZ(1)[1] \bigr) \; .
\]
By \cite[Cor.~3.4.3]{VSF}, for any smooth variety $W$ over $k$, 
the group 
\[
\Hom_{\DeffgM_{\BQ}} \bigl( M(W) , \BZ(1)[1] \bigr)
\]
is naturally isomorphic to the group of global sections
$\Gamma(W, \BG_m)$, tensored with $\BQ$. Therefore,
the inclusion of 
$\Hom_{\DeffgM_{\BQ}} ( M_0(\Xp) \otimes M(Y) , \BZ(1)[1] )$ into
$\Hom_{\DeffgM_{\BQ}} ( M(\Xp \times_k Y) , \BZ(1)[1] )$
corresponds to
\[
\ap^* : \Gamma(\pi_0(\Xp) \times_k Y, \BG_m) \otimes_\BZ \BQ \longto 
\Gamma(\Xp \times_k Y, \BG_m) \otimes_\BZ \BQ \; ,
\]
where $\ap : \Xp \to \pi_0(\Xp)$ is the structure morphism from $\Xp$ to the
scheme $\pi_0(\Xp) := \Spec \Gamma(\Xp, \CO)$ of connected components of $\Xp$. It
is therefore an isomorphism (recall that $\Xp$ is proper).
\end{Proof}

Putting everything together, we thus get the following result.

\begin{Thm} \label{Main}
Assume that all geometric 
irreducible components of $D$ are of genus zero.
Then there is a canonical
morphism
\[
M_1(D) \stackrel{\ip_*}{\longto} M^{!*}_{\ge 2}(\bX) \longonto 
M_2^{!*}(\bX) \; .
\]
\end{Thm}

It will be convenient to interpret this morphism as a one-extension 
$\BE$ in $\DeffgM_{\BQ}$ of the Artin motive
$M_1(D)[-1]$ by $M^{!*}_2(\bX)[-2]$.

\begin{Rem} \label{6E}
(a)~Remark~\ref{6B} shows where to look for a natural candidate for the cone 
of $\BE: M_1(D) \to M^{!*}_2(\bX)$: it should be a canonical
sub-quotient of the motive with compact support $M^c(X)$. \\[0.1cm] 
(b)~Note that the object $M_1(D)$ is trivial 
(and hence so is $\BE$)
if $\Xs$ is smooth. \\[0.1cm]
(c)~Without the assumption on the genus of the
geometric 
irreducible components of $D$, we still get morphisms
\[
M_1(D) \longto M_2^{!*}(\bX) \; ,
\]
by composing $\ip_*: M_1(D) \to M^{!*}_{\ge 1}(\bX)$ with projections $p_2$
from $M^{!*}_{\ge 1}(\bX)$ to its direct factor $M_2^{!*}(\bX)$.
In special cases, the dependence on the choice of the projection $p_2$
may be controlled.
\end{Rem}

\bigskip

%\include{Sec6}
%%%%%%%%%%%%%%%%%%%%%%%%%%%%%%%%%%%%%%%%%%%%%%%%%%%%%%%%%%%%%%%%%%%%%%%
%
%  Section 7
%
%%%%%%%%%%%%%%%%%%%%%%%%%%%%%%%%%%%%%%%%%%%%%%%%%%%%%%%%%%%%%%%%%%%%%%%

\section{Motivic interpretation of
a construction of A.~Caspar}
\label{7}

%%%%%%%%%%%%%%%%%%%%%%%%%%%%%%%%

%%%%%%%%%%%%%%%%%%%%%%%%%%%%%%%%

We keep the geometric situation studied in the previous section:
$\bX$ is a proper
surface over our perfect base field $k$, and 
\[
\vcenter{\xymatrix@R-10pt{
        X \ar@{^{ (}->}[r] \ar@{=}[d] &
        \Xp \ar@{<-^{ )}}[r] \ar[d]_\pi &
        D \ar[d]^\pi \\
        X \ar@{^{ (}->}[r] &
        \Xs \ar@{<-^{ )}}[r] &
        Z
\\}}
\]
is a cartesian diagram
which is a desingularization of the normalization $\Xs$ of $\bX$,
meaning that $\pi$ is proper, $\Xp$ is smooth, $Z$ is finite,
and $D$ a divisor with normal crossings,
whose irreducible components $D_m$ are smooth projective curves. 
Let us start by proving the following result
(compare \cite[Lemma~1.1]{Cs}).

\begin{Lem} \label{7A}
Denote by $\Pic(\Xp)'$ the group of isomorphism classes of line bundles
on $\Xp$, whose restrictions to all $D_m$ are trivial.
Assume that all geometric irreducible
components of $D$ are of genus zero.
Then the map $\jp^* : \Pic(\Xp)' \to \Pic(X)$ induces an isomorphism
\[
\jp^* \otimes \BQ: \Pic(\Xp)' \otimes_\BZ \BQ \isoto 
\Pic(X) \otimes_\BZ \BQ \; .
\]
\end{Lem}

\begin{Proof}
We may assume that our (perfect) base field $k$ is algebraically closed.
Any element in the kernel of $\jp^*: \Pic(\Xp) \to \Pic(X)$ 
is represented by a linear combination $\sum_m a_m D_m$ of the $D_m$.
If the class of $\sum_m a_m D_m$ belongs to $\Pic(\Xp)'$, then its
intersection numbers with all $D_m$ must be zero. Thus the vector
$(a_m)_m$ is in the kernel of the intersection matrix, 
which is invertible (in $\GL_r (\BQ)$) since the intersection pairing on the 
$D_m$ is non-degenerate \cite[p.~6]{M}. Hence $(a_m)_m$ is zero.
For the surjectivity of 
$\jp^* \otimes \BQ$,
observe that $\jp^*: \Pic(\Xp) \to \Pic(X)$ is surjective.
The non-degeneracy of the intersection matrix shows that any
divisor $C$ on $\Xp$ can be modified by a rational linear combination
of the $D_m$ such that the difference $C'$ has trivial intersection numbers
with all the $D_m$. Since these are supposed to be of genus zero,
the restriction of $C'$ to all $D_m$ is principal.
\end{Proof}

\begin{Prop} \label{7a}
Assume that all geometric irreducible
components of $D$ are of genus zero.
There is a canonical morphism
of vector spaces
\[
\Pic(X) \otimes_\BZ \BQ \longto 
\Ext^1_{\DeffgM_{\BQ}} \bigl( M_1(D)[-1] , M^{!*}_0(\bX)(1) \bigr) \; .
\]
Here, $\Ext^1_{\DeffgM_{\BQ}}(\argdot,\argast)$ denotes  
$\Hom_{\DeffgM_{\BQ}}(\argdot,\argast[1])$.
\end{Prop}

\begin{Proof}
As before, denote by $\Pic(\Xp)'$ the group of line bundles
on $\Xp$, whose restrictions to all $D_m$ are trivial. 
Define a morphism
\[
\Pic(\Xp)' \longto 
\Ext^1_{\DeffgM_{\BQ}} \bigl( M_1(D)[-1] , M^{!*}_0(\bX)(1) \bigr) 
\]
by mapping the class of $\CL \in \Pic(\Xp)'$ to the image of 
\[
\BE \in \Ext^1_{\DeffgM_{\BQ}} \bigl( M_1(D)[-1] , M^{!*}_2(\bX)[-2] \bigr) 
\]
(Theorem~\ref{Main})
under $R(c_{\CL}) : M^{!*}_2(\bX)[-2] \to M^{!*}_0(\bX)(1)$
(Variant~\ref{4A'}~(iii)).
Now use Lemma~\ref{7A}.
\end{Proof}

Given a sub-scheme $Z_\infty$ of the finite scheme $Z$, we may consider the
pre-image $D_\infty \subset D$ of $Z_\infty$ under 
$\pi$, and define $M_1(D_\infty)$
as before. It is a direct factor of $M_1(D)$, with a canonical complement. 

\begin{Cor} \label{7b}
Assume that all geometric irreducible
components of $D$ are of genus zero.
There is a canonical morphism of vector spaces
\[
\Pic(X)\otimes_\BZ \BQ \longto 
\Ext^1_{\DeffgM_{\BQ}} \bigl( M_1(D_\infty)[-1] , M^{!*}_0(\bX)(1) \bigr) \; .
\]
\end{Cor}

\begin{Ex} \label{7c}
Here, our base field is equal to $\BQ$.
Let us recall the geo\-metric setting studied in \cite{Cs}.
Let $F$ be a real quadratic extension of $\BQ$ with discriminant $d$. 
Assume that the class number in the narrow sense
of $F$ equals one. Let $X'$ be the \emph{Hilbert modular surface}
of full level associated to $F$ \cite[Sect.~X.4]{vdG}.
Denote by $X^*$ its \emph{Baily--Borel
compactification}, and by $X$ the smooth part of $X'$. All these surfaces
are normal and geometrically connected. The complement of
$X^* - X'$ consists of 
one $\BQ$-rational point, denoted $\infty$ (the \emph{cusp}
of $X^*$). 
The finite sub-scheme $Z := (X^* - X)_{\red}$ includes the cusp, but
also the singularities of $X'$. There is a 
canonical desingularization
\[
\vcenter{\xymatrix@R-10pt{
        X \ar@{^{ (}->}[r] \ar@{=}[d] &
        \Xp \ar@{<-^{ )}}[r] \ar[d]_\pi &
        D \ar[d]^\pi \\
        X \ar@{^{ (}->}[r] &
        \Xs \ar@{<-^{ )}}[r] &
        Z
\\}}
\]
$\Xp$ is a smooth, projective scheme over $\BQ$, and $D$ a divisor with normal cros\-sings, 
whose irreducible components are smooth. Furthermore,
all geometric irreducible components of 
$D$ are of genus zero. The irreducible components
of the pre-image $D_\infty \subset D$ 
of $\infty$ under $\pi$ are isomorphic
to $\BP^1_{\BQ}$, and form a polygon:  
for the complex surface underlying $\Xp$,
this is due to Hirzebruch \cite[Chap.~II]{vdG};
that the statement holds over $\BQ$ follows from 
\cite[Sect.~5]{R}. \\[0.1cm]
(1)~We claim that
the Artin motive $M_1(D_\infty)[-1]$ is canonically isomorphic to
$H_1(D_\infty(\BC),\BZ) \otimes_\BZ \BZ(0)$. (Any of the two orientations of
the polygon $D_\infty$ will thus induce an isomorphism from
$M_1(D_\infty)[-1]$ to $\BZ(0)$.)

Indeed, by the same reasoning
as in Proposition~\ref{6D},
the Artin motive $M_1(D_\infty)[-1]$ equals the kernel of 
\[
\bigoplus_{n < m} M(D_n \cap D_m) \longto \bigoplus_m M_0(D_m) \; ,
\]
where $D_m$ are the components of $D_\infty$. Since $D_\infty$
is a polygon, all $M_0(D_m)$ are equal to $\BZ(0)$,
while the $M_1(D_m)$ are zero.
The $M(D_n \cap D_m)$ are equal to $\BZ(0)$ for consecutive indices $n,m$. 
Hence the kernel in question equals the tensor product of the motive
$\BZ(0)$ with the kernel of the morphism
\[
\bigoplus_{n < m} H_0 \bigl( (D_n \cap D_m)(\BC),\BZ \bigr) \longto 
\bigoplus_m H_0 \bigl( D_m(\BC),\BZ \bigr) 
\]  
of homology groups. \\[0.1cm]
(2)~The variety $\Xp$ being geometrically connected, we have
\[
M^{!*}_0(\bX) = M_0(\Xp) = \BZ(0) \; .
\]
Corollary~\ref{7b} thus yields the following. \\[0.1cm]
(3)~Let $k$ be an extension of $\BQ$. Denote by $X_k$
the base change of $X$ to $k$. Then
there is a canonical morphism $cl_{\KCE}$ mapping
$\Pic(X_k)\otimes_\BZ \BQ$ to
\[
\Ext^1_{\DeffgM_{\BQ}}
              \bigl( H_1(D_\infty(\BC),\BZ) \otimes_\BZ \BZ(0) , \BZ(1) \bigr)
= H^1 \bigl( D_\infty(\BC),k^* \bigr) \otimes_\BZ \BQ \; .
\] 
Any of the two orientations of
the polygon $D_\infty$ thus induces a morphism
\[
cl_{\KCE}: \Pic(X_k)\otimes_\BZ \BQ \longto k^* \otimes_\BZ \BQ \; .
\] 
Indeed, the only point to be verified is the equality
\[
\Ext^1_{\DeffgM_{\BQ}} \bigl( \BZ(0) , \BZ(1) \bigr)
= k^* \otimes_\BZ \BQ \; .
\] 
But this is the content of \cite[Cor.~3.4.3]{VSF}. \\[0.1cm]
(4)~Following the terminology of \cite{Cs},
the image of the class of a line bundle $\CL$ 
under $cl_{\KCE}$ will be called the \emph{Kummer--Chern--Eisenstein extension}
associated to $\CL$. \\[0.1cm]
(5)~Now consider the case $k = F = \BQ (\sqrt{d})$.
Let $\sigma_1,\sigma_2$ be the (real) embeddings of $F$ into $\BC$.
We consider the two line bundles $\CL_i$ 
on $X_F$,
$i = 1,2$, characterized by their factors of automorphy 
``$(\gamma \tau_i + \delta)^2$'' over $\BC$. 
We propose ourselves to identify their images 
under the map $cl_{\KCE}$ from (3). 
To do so, fix an orientation of $D_\infty$. 
Denote by $\varepsilon \in \CO^*_F$
the generator of the totally posi\-tive units. We shall show
(Example~\ref{final}): \emph{if $d$ is a prime congruent to $1$ modulo $4$,
then}
\[
cl_{KCE} (\CL_1 \otimes \CL_2) = 1 \in F^* \otimes_\BZ \BQ
\quad and \quad
cl_{KCE} (\CL_1) = \varepsilon^{\pm 1} \in F^* \otimes_\BZ \BQ  \; .
\]
(The ambiguity concerning the sign in the exponent comes from the choice
of the orientation.) \\[0.1cm]
(6)~This claim implies in particular 
that the realizations of the Kummer--Chern--Eisenstein
extensions $cl_{KCE} (\CL_1)$ and $cl_{KCE} (\CL_2)$ can be identified.
For the $\ell$-adic and Hodge--de Rham realization, this identification
is the content of Caspar's main results \cite[Thm.~2.5, Thm.~3.4]{Cs}. 
Our claim is compatible with \loccit. 
Note that it also implies that the extension 
\[
\BE \in \Ext^1_{\DeffQgM_{\BQ}} \bigl( M_1(D)[-1] , M^{!*}_2(X^*)[-2] \bigr) 
\]
from Theorem~\ref{Main} is non-trivial in the present
geometric situation.
\end{Ex}

In order to prove the claim made in Example~\ref{7c}~(5), let us come back to
the more general situation 
\[
\vcenter{\xymatrix@R-10pt{
        X \ar@{^{ (}->}[r] \ar@{=}[d] &
        \Xp \ar@{<-^{ )}}[r]^{\ip} \ar[d]_\pi &
        D \ar[d]^\pi \\
        X \ar@{^{ (}->}[r] &
        \Xs \ar@{<-^{ )}}[r] &
        Z
\\}}
\]
considered in the beginning of this section.
In particular, the irreducible components
$D_m$ of $D$ are supposed smooth (and projective),
but not ne\-cessarily of genus zero.
We need to generalize the construction of the cup product
with the first Chern class of a line bundle. Recall that 
for a smooth and projective variety $Y$,
the vector space $CH^1(Y) = \Pic(Y) \otimes_\BZ \BQ$ equals 
\[
\Hom_{CHM(k)_{\BQ}} \bigl( \BL , h(Y) \bigr) = 
\Hom_{\DeffgM_{\BQ}} \bigl( M(Y) , \BZ(1)[2] \bigr) \; .
\]
In fact, Voevodsky \cite[Cor.~3.4.3]{VSF}
proved the following result.

\begin{Thm} \label{7d}
Let $Y \in Sm/k$. For any $j \in \BZ$, 
there is a canoni\-cal isomorphism
\[
H_{Zar}^{j-1}(Y,\BG_m) \isoto
\Hom_{\DeffgM} \bigl( M(Y) , \BZ(1)[j] \bigr) \; ,
\]
which is contravariantly functorial in $Y$.
\end{Thm}

In particular, we then have $\Pic(Y) =
\Hom_{\DeffgM} \bigl( M(Y) , \BZ(1)[2] \bigr)$. 
It follows from the construction of \loccit \
that for $Y$ smooth and projective, the tensor product of
this isomorphism with $\BQ$ is the one we used in Section~\ref{4} to produce
morphisms $\BL \to h(Y)$ of Chow motives. 
Analyzing more closely the ingredients of Voevodsky's proof,
we are able to show the following.

\begin{Prop} \label{7e}
(i)~There is a canonical isomorphism 
\[
\Pic(D) \isoto \Hom_{\DeffgM} \bigl( M(D) , \BZ(1)[2] \bigr) \; .
\]
(ii)~The diagram
\[
\vcenter{\xymatrix@R-10pt{
        \Pic(D) \ar[r]^-{\cong} &
        \Hom_{\DeffgM} \bigl( M(D) , \BZ(1)[2] \bigr) \\
        \Pic(\Xp) \ar[r]^-{\cong} 
                                        \ar[u]^{{\ip}^*} &
        \Hom_{\DeffgM} \bigl( M(\Xp) , \BZ(1)[2] \bigr) 
                                        \ar[u]_{{\ip}^*}
\\}}
\]
commutes. \\[0.1cm]
(iii)~Denote by $\ip_m$ the inclusion of $D_m$ into $D$.
Then for all $m$, the diagram
\[
\vcenter{\xymatrix@R-10pt{
        \Pic(D_m) \ar[r]^-{\cong} &
        \Hom_{\DeffgM} \bigl( M(D_m) , \BZ(1)[2] \bigr) \\
        \Pic(D) \ar[r]^-{\cong} 
                                        \ar[u]^{{\ip_m}^*} &
        \Hom_{\DeffgM} \bigl( M(D) , \BZ(1)[2] \bigr) 
                                        \ar[u]_{{\ip_m}^*}
\\}}
\]
commutes.
\end{Prop}

\begin{Proof}
Recall (see the introduction to Section~\ref{5})
that $M = \RC \circ L$,
and that $\RC: D^-(\ShN) \to \DM$ is left adjoint to the inclusion
of $\DM$ into $D^-(\ShN)$.
It follows that for any Nisnevich sheaf with transfers $G$,
any integer $r$, and any $Y \in Sch/k$, we have
\[
\Hom_{\DM} \bigl( M(Y) , G[r] \bigr) = 
\Hom_{D^-(\ShN)} \bigl( L(Y) , G[r] \bigr) \; .
\]
Note that if $Y$ is smooth, then $L(Y)$ is the Nisnevich sheaf 
with transfers represented by $Y$,
hence by Yoneda's Lemma, 
\[
\Hom_{\ShN} \bigl( L(Y) , G \bigr) = \Gamma(Y,G) \; .
\] 
By definition of $L$,
the sequence
\[
0 \longto \bigoplus_{n < m} L(D_n \cap D_m) \longto \bigoplus_n L(D_n) 
\longto L(D) \longto 0 
\]
is exact (even as a sequence of presheaves --- recall that the $D_n$ are 
the irreducible components of $D$).
This shows that 
\[
\Hom_{\ShN} \bigl( L(D) , G \bigr) =
\ker \bigl( \prod_n \Gamma(D_n,G) \longto 
\prod_{n < m} \Gamma(D_n \cap D_m,G) \bigr) \; .
\]
For any open subset $U$ of $D$, the formula
\[
\Gamma(U,\FH^0(G)) := 
\ker \bigl( \prod_n \Gamma(D_n \cap U,G) \longto 
\prod_{n < m} \Gamma(D_n \cap D_m \cap U,G) \bigr)
\]
\emph{defines} a functor on $\ShN$.
Letting $U$ vary, we get a left exact functor
\[
\FH^0 : \ShN \longto \ShDZ \; ,
\] 
where we denote by $\ShDZ$ the category of Zariski sheaves 
with values in Abelian groups on the topological 
space underlying $D$. 
We claim that there are natural morphisms
\[
H^r_{Zar} \bigl( D,\FH^0(G) \bigr) \longto 
\Hom_{D^-(\ShN)} \bigl( L(D) , G[r] \bigr) 
\]
for any Nisnevich sheaf with transfers $G$.
Observe that by what was said before, this is a natural 
isomorphism for $r = 0$.
The morphisms in question will be defined as the boundaries in 
a spectral sequence
\[
H^p_{Zar} \bigl( D,R^q ( \FH^0 )(G) \bigr)
\Longrightarrow 
\Hom_{D^-(\ShN)} \bigl( L(D) , G[p+q] \bigr) 
\]
which we construct now.
The category $\ShN$ has sufficiently 
many injectives \cite[Lemma~3.1.7]{VSF}.
Hence the existence of the spectral sequence is equivalent to  
\[
\quad H^r_{Zar} \bigl( D,\FH^0(I) \bigr) = 0 \; , \; r \ge 1 \; ,
\]
for any injective $I \in \ShN$.
The proof of this vanishing
is a faithful imitation of the proof of 
\cite[Prop.~3.1.8]{VSF}; note that the vital
ingredient of \loccit \ is \cite[Prop.~3.1.3]{VSF},
which is valid without any smoothness assumptions.

Let us now specialize to the case $G = \BG_m$ and $r = 1$.
For two indices $n < m$, denote
by $\ip_{n,m}$ the inclusion of $D_n \cap D_m$ into $D$. 
The short exact sequence of Zariski sheaves on $D$
\[
(\ast) \quad \quad
1 \longto \BG_{m,D} \longto \prod_n \ip_{n,*} \BG_{m,D_n}
\longto \prod_{n < m} \ip_{n,m,*} \BG_{m,D_n \cap D_m} \longto 1  
\]
shows that $\BG_{m,D} = \FH^0(\BG_m)$.
Hence the above construction yields 
\[
\Pic(D) = H^1_{Zar} \bigl( D,\BG_m \bigr) \longto 
\Hom_{D^-(\ShN)} \bigl( L(D) , \BG_m[1] \bigr) \; .
\]
But by \cite[Thm.~3.4.2]{VSF}, there is a canonical isomorphism
$\BZ(1)[1] \cong \BG_m$ in $\DM \subset D^-(\ShN)$.
Altogether, we get the required morphism
\[
\Pic(D) \longto 
\Hom_{\DeffgM} \bigl( M(D) , \BZ(1)[2] \bigr) \; .
\]
By construction, it is compatible with the isomorphisms from Theorem~\ref{7d}
(for $j = 2$)
under morphisms of schemes $Y \to D$ and $D \to Y$, for $Y \in Sm/k$.

It remains to show that $\Pic(D) \to \Hom_{\DeffgM} ( M(D) , \BZ(1)[2])$
is in fact an isomorphism. But this follows easily from the Five Lemma,
from the long exact Zariski cohomology sequence induced by $(\ast)$,
and the long exact $\Hom_{\DeffgM} ( \bullet , \BZ(1)[1])$-sequence
induced by the exact triangle 
\[
\bigoplus_{n < m} M(D_n \cap D_m) \longto \bigoplus_n M(D_n) 
\longto M(D) \longto \bigoplus_{n < m} M(D_n \cap D_m)[1] \; ,
\]
and from Theorem~\ref{7d}.
\end{Proof}

\forget{
The short exact sequence induces a long exact
cohomology sequence, which shows that $\Pic(D)'$ equals the cokernel
of the map
\[
\prod_m H^0 \bigl( D_m , \CO_{D_m}^* \bigr) \longto 
\prod_{n < m} H^0 \bigl( D_{n,m} , \CO_{D_{n,m}}^* \bigr) \; .
\]
The exact triangle induces an exact 
$\Hom_{\DeffgM} ( \bullet , \BZ(1)[1] )$-triangle, which shows that
$\Hom_{\DeffgM} \bigl( M(D) , \BZ(1)[2] \bigr)'$ equals the cokernel
of the map
\[
\prod_m \Hom_{\DeffgM} \bigl( M(D_m) , \BZ(1)[1] \bigr) \longto
\prod_{n < m} \Hom_{\DeffgM} \bigl( M(D_{n,m}) , \BZ(1)[1] \bigr) \; .
\]
Now use Theorem~\ref{7d}.
}

\begin{Rem}
We leave it to the reader to prove that the conclusions of Proposition~\ref{7e}
are in fact true whenever $D$ is a normal crossing divisor in $\Xp \in Sm / k$,
with smooth irreducible components $D_m$.
\end{Rem}

\forget{
\begin{Cor} \label{7f}
Denote by $\Pic(D)'$ and $\Pic(\Xp)'$
the groups of classes of line bundles
on $D$ and $\Xp$ respectively, whose restrictions to all $D_m$ are trivial.
Similarly, denote by
\[
\Hom_{\DeffgM} \bigl( M(D) , \BZ(1)[2] \bigr)' \quad \text{and} \quad
\Hom_{\DeffgM} \bigl( M(\Xp) , \BZ(1)[2] \bigr)'
\]
the groups of morphisms $M(D) \to \BZ(1)[2]$ 
and $M(\Xp) \to \BZ(1)[2]$ respectively, whose restrictions to
all $M(D_m)$ are trivial. \\[0.1cm]
(i)~There is a canonical isomorphism 
\[
\Pic(D)' \isoto \Hom_{\DeffgM} \bigl( M(D) , \BZ(1)[2] \bigr)' \; .
\]
(ii)~The diagram
\[
\vcenter{\xymatrix@R-10pt{
        \Pic(D)' \ar[r]^-{\cong} &
        \Hom_{\DeffgM} \bigl( M(D) , \BZ(1)[2] \bigr)' \\
        \Pic(\Xp)' \ar[r]^-{\cong} 
                                        \ar[u]^{{\ip}^*} &
        \Hom_{\DeffgM} \bigl( M(\Xp) , \BZ(1)[2] \bigr)' 
                                        \ar[u]_{{\ip}^*}
\\}}
\]
commutes.  
\end{Cor}
}

For any line bundle $\CK$ on $D$, we can now define a morphism
\[
R(c_\CK) : M(D) \longto M(D)(1)[2]
\]
in complete analogy to the smooth projective case, namely as the composition
\[
M(D) \stackrel{\Delta_*}{\longto} M(D) \otimes M(D) 
\stackrel{\id_{D,*} \otimes [\CK]}{\longto} M(D)(1)[2] 
\]
($\Delta :=$ the diagonal embedding $D \into D \times_k D$).

\begin{Cor} \label{7g}
(i)~Let $\CL$ be a line bundle on $\Xp$. Then the diagram
\[
\vcenter{\xymatrix@R-10pt{
        M(D) \ar[r]^-{R(c_{\ip^* \! \CL})} \ar[d]_{\ip_*} &
        M(D)(1)[2] \ar[d]^{\ip_*(1)[2]} \\
        M(\Xp) \ar[r]^-{R(c_{\CL})} &
        M(\Xp)(1)[2]
\\}}
\]
commutes. \\[0.1cm]
(ii)~Let $\CK$ be a line bundle on $D$.
Then for all $m$, the diagram
\[
\vcenter{\xymatrix@R-10pt{
        M(D_m) \ar[r]^-{R(c_{\ip_m^* \! \CK})} \ar[d]_{\ip_{m,*}} &
        M(D_m)(1)[2] \ar[d]^{\ip_{m,*}(1)[2]} \\
        M(D) \ar[r]^-{R(c_{\CK})} &
        M(D)(1)[2]
\\}}
\]
commutes.
\end{Cor}

\begin{Cor} \label{7h}
Let $\CK$ be a line bundle on $D$, whose restrictions to all $D_m$
are trivial. Then $R(c_{\CK}):M(D) \to M(D)(1)[2]$
factors uniquely through a morphism $R(c_{\CK}): M_{\le 1}(D) \to M(D)(1)[2]$.
\end{Cor}

\begin{Proof}
Recall that $M_{\le 1}(D)$ is the categorial quotient of $M(D)$
by $M_2(D)$. Our claim thus
follows from Corollary~\ref{7g}~(ii), Proposition~\ref{7e}~(iii)
and the equation $M_2(D) = \oplus_m M_2(D_m)$.
\end{Proof}

Composition with the monomorphism $M_1(D) \into M_{\le 1}(D)$
and the epimorphism $M(D)(1)[2] \onto M_0(D)(1)[2]$ thus yields a map
\[
cl_D: \Pic(D)' \otimes_{\BZ} \BQ \longto 
\Ext^1_{\DeffgM_{\BQ}} \bigl( M_1(D)[-1] , M_0(D)(1) \bigr) \; .
\]

\begin{Prop} \label{7i}
Assume that all geometric irreducible
components of $D$ are of genus zero.
Then the morphism
\[
cl_X: \Pic(X) \otimes_\BZ \BQ \longto 
\Ext^1_{\DeffgM_{\BQ}} \bigl( M_1(D)[-1] , M^{!*}_0(\bX)(1) \bigr) 
\]
of Proposition~\ref{7a} factors canonically through $cl_D$.
More precisely, the diagram
\[
\vcenter{\xymatrix@R-10pt{
        \Pic(D)' \otimes_\BZ \BQ \ar[rr]^-{cl_D} &&
        \Ext^1 \bigl( M_1(D)[-1] , M_0(D)(1) \bigr) 
            \ar[d]^{{\ip}_*} \\
        \Pic(\Xp)' \otimes_\BZ \BQ \ar[r]^-{\cong}_-{\ref{7A}} 
            \ar[u]^{{\ip}^*} &
        \Pic(X) \otimes_\BZ \BQ \ar[r]^-{cl_X} &
        \Ext^1 \bigl( M_1(D)[-1] , M^{!*}_0(\bX)(1) \bigr) 
                                        \\}}
\]
commutes, where we abbreviated $\Ext^1 := \Ext^1_{\DeffgM_{\BQ}}$.
\end{Prop}

\begin{Proof}
Let $\CL$ be a line bundle on $X$.
Recall that the morphism of Proposition~\ref{7a} maps  
the class of $\CL$ to
the image of 
\[
\BE \in \Ext^1_{\DeffgM_{\BQ}} \bigl( M_1(D)[-1] , M^{!*}_2(\bX)[-2] \bigr) 
\]
(Theorem~\ref{Main})
under $R(c_{\CL}) : M^{!*}_2(\bX)[-2] \to M^{!*}_0(\bX)(1)$
(Variant~\ref{4A'}~(iii)), where by abuse of notation we denote by $\CL$ 
also the unique extension of $\CL$ to $\Pic(\Xp)' \otimes_\BZ \BQ$
(Lemma~\ref{7A}). Our claim thus follows from Corollary~\ref{7g}~(i).
\end{Proof}

\begin{Ex} \label{final}
Let us reconsider the situation from Example~\ref{7c}, and 
prove the claim made in \ref{7c}~(5). The polygon $D_\infty$ is geometrically
connected, therefore $M_0(D_\infty) \to M^{!*}_0(\bX)$ 
is an isomorphism (both sides equal $\BZ(0)$). 
By Proposition~\ref{7i}, the morphism
\[
cl_{\KCE}: \Pic(X_k) \otimes_\BZ \BQ \longto
H^1 \bigl( D_\infty(\BC),k^* \bigr) \otimes_\BZ \BQ 
\] 
factors through $cl_{D_\infty}$, where
\[
cl_{D_\infty}: \Pic(D_{\infty,k})' \otimes_\BZ \BQ \longto
H^1 \bigl( D_\infty(\BC),k^* \bigr) \otimes_\BZ \BQ \; .
\] 
Using the long exact Zariski cohomology sequence induced by 
\[
1 \longto \BG_{m,D_\infty} \longto \prod_n \ip_{n,*} \BG_{m,D_n}
\longto \prod_{n < m} \ip_{n,m,*} \BG_{m,D_n \cap D_m} \longto 1 
\]
and the calculation of \ref{7c}~(1),
one sees that $cl_{D_\infty}$ is in fact an isomorphism.
Any of the two orientations of
the polygon $D_\infty$ thus induces an isomorphism
\[
cl_{D_\infty}: \Pic(D_{\infty,k})' \otimes_\BZ \BQ \isoto k^* \otimes_\BZ \BQ \; .
\]  
Checking the definitions, we can identify $cl_{D_\infty}$: we fix a 
point $x_0 \in D_\infty(k)$. It lies on a component $D_{m_0}$.
For any line bundle $\CK$ on $D_{\infty,k}$ with trivial restrictions
to all $D_{m,k}$, we fix an element $s$ in the fibre $\CK_{x_0}$.
The restriction $\Gamma(D_{m_0,k},\CK) \to \CK_{x_0}$ being an isomorphism,
$s$ can be uniquely extended to the whole of $D_{m_0,k}$. We restrict this
extension to the ($k$-rational)
point $x_1$ which is the intersection of $D_{m_0}$ with the ``next''
component (in the sense of the chosen orientation). We repeat the
process until we are again on $D_{m_0}$. Restriction to $\CK_{x_0}$
gives a non-zero multiple $c \cdot s$, and we have $cl_{D_\infty}([\CK]) = c$.
   
In order to prove the claim made in \ref{7c}~(5), one needs to apply
this recipe to the line bundles $\CK_i$ obtained by restricting 
to $D_{\infty,F}$ the unique extensions of $\CL_i$ to 
$\Pic(\Xp_F)' \otimes_\BZ \BQ$, $i = 1,2$. But this is exactly the content
of \cite[Lemma~1.2]{Cs}.
\end{Ex}

\bigskip

\end{document}